\documentclass[reqno, 11pt]{amsart}
\usepackage{amssymb, amsthm, amsmath}
\numberwithin{equation}{section}

\newcommand{\K}{\mathfrak{K}} 
\renewcommand{\k}{\mathfrak{k}} 
\newcommand{\R}{\mathbb{R}} 

\newcommand{\Z}{\mathbb{Z}} 
\newcommand{\N}{\mathbb{N}}
 
\newcommand{\Q}{\mathbb{Q}}
\newcommand{\x}{\ma{x}} 
\newcommand{\X}{\ma{X}}

\newcommand{\rom}{\mathrm} 
\newcommand{\bfP}{\mathbb{P}}
 
\newcommand{\ma}{\mathbf}
\newcommand{\ben}{\begin{enumerate}}
\newcommand{\een}{\end{enumerate}} 
\newcommand{\eit}{\begin{itemize}}
\newcommand{\beq}{\begin{equation}} 
\newcommand{\eeq}{\end{equation}}
 
\newcommand{\ve}{\varepsilon}
 
\newcommand{\mcal}{\mathcal}
\newcommand{\lab}{\label} 
\newcommand{\al}{\alpha}
\newcommand{\D}{\Delta_Q} 
\newcommand{\del}{\delta}

\newcommand{\be}{\beta}

\newcommand{\h}{\mathrm{gcd}}
\newtheorem{theorem}{Theorem} 
\newtheorem{lemma}{Lemma}
\newtheorem{pro}{Proposition} 
\newtheorem{cor}{Corollary}
 
\renewcommand{\mod}{\hspace{-0.25cm}\pmod}
\newcommand{\tmod}{\hspace{-0.1cm}\pmod}
\newcommand{\colt}[2]{\genfrac{}{}{0pt}{1}{#1}{#2}}

\newcommand{\hcf}{\mathrm{gcd}}
\renewcommand{\c}{\ma{c}} 
 
\newcommand{\y}{\ma{y}}
\newcommand{\z}{\ma{z}} 
\renewcommand{\u}{\ma{u}} 

\renewcommand{\leq}{\leqslant} 
\renewcommand{\le}{\leqslant} 
\renewcommand{\geq}{\geqslant}
\renewcommand{\ge}{\geqslant}
\renewcommand{\d}{\mathrm{d}}
\newcommand{\Mq}{\|Q\|} 
\newcommand{\mq}{m(Q)}
\newcommand{\CC}{\mathcal{C}} 
\renewcommand{\ss}{\mathfrak{S}}

\DeclareMathOperator{\supp}{supp}
\theoremstyle{definition} 
\newtheorem*{ack}{Acknowledgement}
\newtheorem*{notation}{Notation}

\begin{document}

\title{On the Representation of integers by quadratic forms}

\author{T.D. Browning}
\author{R. Dietmann}

\address{Institut f\"ur Algebra und Zahlentheorie,
Lehrstuhl f\"ur Zahlentheorie,
Pfaffenwaldring 57,
D-70569 Stuttgart}
\email{dietmarr@mathematik.uni-stuttgart.de}

\address{School of Mathematics,  
University of Bristol, Bristol BS8 1TW}
\email{t.d.browning@bristol.ac.uk}

%\date{\today}

\begin{abstract}
Let $n \geq 4$, and let  $Q\in
\Z[X_1,\ldots,X_n]$ be a  non-singular quadratic form.
When $Q$ is indefinite we provide new upper bounds for the 
least non-trivial integral solution to the equation $Q=0$, and when $Q$
is positive definite we provide improved upper bounds for the least positive
integer $k$ for which the equation $Q=k$ is insoluble
in integers, despite being soluble modulo every prime power.
\end{abstract}

\maketitle

\section{Introduction}\lab{intro}

Let $n \geq 3$ and let $Q\in \Z[X_1,\ldots,X_n]$  
be a non-singular quadratic form, with discriminant $\D$.  
Given a non-negative integer $k$, this paper 
is concerned with the locus of points
$$
\mcal{S}(k;Q):=\{\x=(x_1,\ldots,x_n)\in \Z^n: ~\x\neq\ma{0},~Q(\x)=k\}.  
$$
Our basic motivation is the problem of deciding whether or not $\mcal{S}(k;Q)$ is empty, for given $k$ and
$Q$.  We shall address this question in two rather different
contexts: firstly when $Q$ is indefinite and $k=0$, and secondly,
when $k$ is large and $Q$ is positive definite.  
Throughout this paper we shall work 
with classically integral quadratic forms $Q$. 
Thus there is a symmetric matrix $\ma{A}=(A_{ij})_{1\leq i,j\leq n}$,
with coefficients in $\Z$, such that
$$
Q(\X)=\X^T \ma{A}\X
$$
and $\D= \det \ma{A}$.  We shall write
$$
\Mq:=\max_{1\leq i,j \leq n}|A_{ij}|
$$
for the height of the quadratic form $Q$.  
We are now ready to discuss the decidability question for
$\mcal{S}(k;Q)$, for which we distinguish between the indefinite and
positive definite cases.

\subsection{Indefinite forms}

In this section we consider the problem of deciding when a given indefinite quadratic
form represents zero. 
Now it is well-known that the Hasse--Minkowski theorem gives a 
finite procedure for determining whether or not the set 
$$
\mcal{S}(Q):=\mcal{S}(0;Q)
$$ 
is empty.  An alternative procedure arises from 
providing an effective upper bound for the maximum modulus of the smallest
element of $\mcal{S}(Q)$, in terms of the coefficients of $Q$.

Let $\Lambda_n(Q)$ be minimal with the property that when
$\mcal{S}(Q)$ is non-empty, it contains an element with $
|\x|  \leq \Lambda_n(Q)$. Here, as throughout our work, we 
take $|\z|$ to be the norm $\max_{1\leq i \leq n}|z_i|$, for
any $\z\in \R^n$.  In this context there is a rather old result due to Cassels
\cite{cassels}, which shows that  
\beq\lab{55}
\Lambda_n(Q) \leq c_n\Mq^{(n-1)/2},
\eeq
with $c_n=(\frac 1{2}(3n^2+n-10)(n-1)!^2)^{(n-1)/2}.$
A simpler proof of \eqref{55} has been given by Davenport
\cite{dav-hermite}, with the sharper value 
$(\sqrt{2}n\gamma_{n-1})^{(n-1)/2}$
for the constant $c_n$.  Here, $\gamma_m$ is 
Hermite's constant, defined as the upper bound of the
minima of all positive definite quadratic forms in $m$ variables, of
determinant $1$. 
While the precise value of the constant $c_n$ in
\eqref{55} is perhaps unimportant, the exponent of $\Mq$ has much
more significance.  In fact it was  shown to be best
possible by Kneser \cite{cassels'}, via the ingenious example
$$
Q_0(\X)=X_1^2-(X_2-cX_1)^2-\cdots-(X_n-cX_{n-1})^2,
$$
for any integer $c\geq 3$. It is self-evident 
that $Q_0$ is a non-singular indefinite
quadratic form, with height $\|Q_0\|=c^2$.  Moreover we clearly have
$Q_0(\ma{a})=0$, where
$$
\ma{a}=(1, c-1,c^2-c,\ldots,c^{n-1}-c^{n-2}).
$$
A little thought reveals that $\ma{a}$ is the unique solution
to the equation $Q_0=0$, with least norm and positive first component.
On noting that 
%Here any non-zero integer solution
%$\ma{a}\in \Z^n$ to the equation $Q_0=0$ must have $a_1\neq 0$, and 
%$$
%a_n=b_n+cb_{n-1}+\cdots+ c^{n-2}b_2+c^{n-1}a_1,
%$$
%where $b_j=a_j-ca_{j-1}$ for $j \geq 2$, and $b_n^2+\cdots
%+b_2^2=a_1^2$.  In particular the smallest absolute value of $a_n$ is
$c^{n-1}-c^{n-2}>\frac{1}{2}c^{n-1}=\frac{1}{2}\|Q_0\|^{(n-1)/2}$, since
$c\geq 3$, we therefore conclude that the exponent of $\Mq$ in
\eqref{55} is, in general, best possible.

The first goal of this paper is to show that Cassels'
exponent can be sharpened substantially 
when suitable assumptions are made about the form $Q$.
That such improvements are available already follows from the work of 
Schlickewei \cite{schl}.  Given a positive integer $d\leq n/2$, it is shown
in this work that there is a constant $\tilde{c}_{n}>0$ depending only on $n$, such that 
$$
\Lambda_n(Q)\leq \tilde{c}_{n} \Mq^{(n-d)/(2d)},
$$
provided that $Q$ vanishes on a $d$-dimensional subspace
of $\Q^n$.  One retrieves \eqref{55} by taking $d=1$ in
this inequality. In our present work we shall derive alternative hypotheses under which
the exponent of Cassels can be beaten.
Suppose for the moment that $Q$ is diagonal, with $\D\neq 0$, so that
\beq\lab{Q} 
Q(\X)=A_1X_1^2+\cdots+ A_nX_n^2, 
\eeq 
for non-zero $A_1,\ldots,A_n\in \Z$ not all of the same sign.  In
particular we now have $\D=A_1\cdots A_n$ and $\Mq=\max_{1\leq i \leq n}|A_{i}|$.
In this setting it is worth drawing attention to the work of
Ou and Williams \cite{ou}, who have shown that
if the equation $Q=0$ is non-trivially soluble in integers $x_1, \ldots, x_n$, then there is 
a solution satisfying 
$$
|A_1|x_1^2 + \cdots + |A_n|x_n^2 \leq 2|\D|. 
$$
In view of the fact that $|A_i|^{-1}|\D|\leq \Mq^{n-1}$,
for $1 \leq i \leq n$, this result implies that 
\beq\lab{o-w}
\Lambda_n(Q) \leq \sqrt{2}\Mq^{(n-1)/2},
\eeq  
for any indefinite quadratic form of the shape \eqref{Q}.
Thus not only do we get an improvement over \eqref{55} in the value of the
constant, but we also get an improvement over the term $\Mq^{(n-1)/2}$,
if the coefficients of $Q$ don't all have the same order of magnitude.
It seems to be an open question whether or not there exists a version 
of \eqref{55}, for arbitrary indefinite quadratic forms, in which the
constant $c_n$ is actually independent of $n$.

One obvious way of beating Cassels' bound for large values of $n$
arises through setting $n-5$ of the variables equal to zero, and
then applying \eqref{55} to the resulting form.  In view of Meyer's theorem
this form is guaranteed to have at least one non-trivial solution
provided that it is indefinite, so that Cassels' inequality is applicable.
Note that a quinary quadratic form of rank at most $4$ trivially
represents zero. That the procedure of eliminating $n-5$ of the variables
is not always available, is demonstrated by Kneser's example above. The
form $Q_0(X_1,\ldots,X_{j-1},0,X_{j+1},\ldots,X_n)$ is negative
definite for any choice of $1 \leq j \leq n$.  
It is not hard to see, however, that this procedure will always be
successful for diagonal indefinite quadratic forms. 
Thus, given any quadratic form of the shape \eqref{Q}, with $n \geq 5$
and $A_1,\ldots,A_n \in \Z$ non-zero and not all of the same sign, it is
always possible to set $n-5$ of the variables equal to zero in such a
way that the resulting quinary form is indefinite. This rather simple
observation, that was drawn to the authors' attention by Professor Heath-Brown, 
leads to the following improvement of \eqref{55} for diagonal
quadratic forms in at least five variables.

\begin{theorem}\lab{diag}
Let $n\geq 5$ and assume that $Q \in \Z[X_1,\ldots,X_n]$ is a
diagonal indefinite quadratic form. Then we have
$$
\Lambda_n(Q) \leq \sqrt{2}  \Mq^2.
$$
\end{theorem}

Here, rather than applying \eqref{55} in the case $n=5$, we have instead
employed the inequality \eqref{o-w} of Ou and Williams to get an even 
sharper result.  
Returning to the general case, let $n \geq 5$ and let $Q \in \Z[X_1,\ldots,X_n]$ be a
non-singular indefinite quadratic form, with underlying matrix $\ma{A}$.  
There does not appear to be a very clean description of the general class $\mcal{Q}(n)$, say,  of
all quadratic forms for which one may set $n-5$ of the variables equal to zero
in such a way that the resulting quinary form is indefinite.
This is unfortunate, since it is clear from our argument above that
$\Lambda_n(Q)\ll \Mq^{2}$, for any $Q \in \mcal{Q}(n)$.  
We shall say no more about the class $\mcal{Q}(n)$ here, save to
observe that a necessary condition for the form $Q$ to be contained in $\mcal{Q}(n)$ is that 
the set of all principal minors of $\mathbf{A}$, with $\D=\det \ma{A}$ removed,
should contain elements of opposite sign.  
Instead, we shall adopt a rather different approach.

Let us define $\Lambda_n^\dagger(Q)$ to have the property that
when there is a vector $\x\in \Z^n$ for which
$Q(\x)=0$ and $x_1\neq 0$, then there exists such a vector with 
$|\x|  \leq \Lambda_n^\dagger(Q)$.  In particular it is plain that
we always have $\Lambda_n(Q)\leq \Lambda_n^\dagger(Q)$.
Let $\lambda_1,\ldots,\lambda_n \in \R$ denote the eigenvalues of $\ma{A}$, and define
\beq\lab{min-Q}
\mq:=\min_{1\leq i\leq n}|\lambda_i|
\eeq
to be the minimum modulus of these eigenvalues.
We shall occasionally appeal to the well-known
inequality
\beq\lab{eigen}
\max_{1 \leq i \leq n}|\lambda_i| \leq n \|Q\|,
\eeq
from which it follows that $\mq\leq n \Mq$.
The following result is our main contribution to the theory of
indefinite quadratic forms.

\begin{theorem}\lab{c}
Let $n\geq 5$ and assume that $Q \in \Z[X_1,\ldots,X_n]$ is an
indefinite quadratic form, with discriminant
$\D\neq 0$.  Then we have
$$
\Lambda_n^\dagger(Q) \ll_\ve
\mq^{-1/2}\big(|\D|^{1+2\be_Q} \Mq^{n+\ve}\big)^{1/(n-3-\al_n)},
$$
for any $\ve>0$, where 
\beq\lab{parity}
\al_n:=\left\{
\begin{array}{ll}
1, & \mbox{if $n$ is even},\\
0, & \mbox{if $n$ is odd,}
\end{array}
\right.
\eeq 
and
\beq\lab{8-q}
\be_Q:=\left\{
\begin{array}{ll}
0, & \mbox{if $\D$ is odd and square-free},\\
1/(n-4), & \mbox{otherwise.}
\end{array}
\right.
\eeq 
\end{theorem}

All of the implied constants that appear in our
work are effectively computable.  They will be permitted to depend at most upon
$n$.  Any further dependence will
be explicitly indicated by appropriate parameters in subscript.
As will become apparent at the close of \S \ref{2}, one may actually take $\be_Q=0$
in the theorem for a rather less restrictive class of quadratic
forms.

It is interesting to place Theorem \ref{c} in the context of our
discussion of Kneser's quadratic form $Q_0$.
This form has discriminant of modulus $1$, and height $\|Q_0\|=c^2$.
Moreover, it turns out that $m(Q_0)\gg c^{2-2n}$, in the notation of
\eqref{min-Q}.  To see this we follow an argument suggested to us by
Professor Heath-Brown.  It begins with the observation that $Q_0$ has underlying
matrix $\ma{A}=\ma{B}^T\ma{D}\ma{B}$, where 
$\ma{D}=\mathrm{Diag}(1,-1,\ldots,-1)$, $\ma{B}=\ma{I}-\ma{C}$, and 
$\ma{C}$ is the matrix with $c$'s just below the diagonal and zero
everywhere else.  But then $\ma{B}^{-1}=\ma{I}+\ma{C}+\cdots
+\ma{C}^{n-1}$, and it is not hard to see that the maximum modulus of
any element of the matrix $\ma{A}^{-1}=\ma{B}^{-1}\ma{D}\ma{B}^{-T}$
is $O(c^{2n-2})$.  Thus, if $\lambda_1,\ldots,\lambda_n$ are the eigenvalues
of $\ma{A}$, it follows from \eqref{eigen} that $m(Q_0)^{-1}\leq \max_{1\leq i \leq
  n}|\lambda_i^{-1}| \ll c^{2n-2}$. This establishes the claim that
$m(Q_0)\gg c^{2-2n}$, and we may now put all of this together to deduce from Theorem \ref{c}
that 
$$
\Lambda_n(Q_0)\ll_\ve
c^{n-1}c^{2n/(n-3-\al_n)+\ve}\leq \|Q_0\|^{(n-1)/2+n/(n-4)+\ve}.
$$
This is, as was to be expected, weaker than the inequality implied by \eqref{55}.

At first glance, it is perhaps not obvious that Theorem \ref{c} ever
improves upon Cassels' result.  To see that it does, it will be
convenient to derive a weaker version of Theorem \ref{c},
in which the term $\mq$ does not appear. It follows from \eqref{eigen}
that $\mq=\mq|\D|/|\lambda_1\cdots \lambda_n|\gg |\D| \Mq^{1-n}$.
Once inserted into Theorem \ref{c}, this yields the following result.

\begin{cor}\lab{cor2}
Let $n\geq 5$ and assume that $Q \in \Z[X_1,\ldots,X_n]$ is an
indefinite quadratic form, with discriminant $\D\neq 0$.   Then we have
$$
\Lambda_n^\dagger(Q) \ll_\ve
\Mq^{(n-1)/2}\Big(\frac{\Mq^{n+\ve}}{|\D|^{(n-5-\al_n-4\be_Q)/2}}\Big)^{1/(n-3-\al_n)}, 
$$
for any $\ve>0$, 
where $\al_n$ is given by \eqref{parity} and $\be_Q$ is given by \eqref{8-q}.
\end{cor}

Suppose for the moment that $\theta\in \R$ is chosen so that
$|\D|=\Mq^{\theta}$, and that $\D$ is odd and square-free.  
Then it follows from Corollary
\ref{cor2} that we get a saving over Cassels' bound for $\Lambda_n(Q)$
as soon as 
$$
\theta>\frac{2n}{n-5-\al_n}.
$$
In particular we must have $n \geq 7$, since $0\leq \theta \leq n$. 
As is well-known, the discriminant of a generic quadratic form has
the same order of magnitude as the 
$n$th power of its height. Thus, for typical indefinite quadratic
forms, one should be able to take $\theta=n$ in 
the above analysis. In  this favourable setting we get a
decent saving over the bound of Cassels as soon as $n\geq 9$.

Whereas we have so far only been interested in Theorem \ref{c} on
the grounds that $\Lambda_n(Q)\leq \Lambda_n^\dagger(Q)$, the 
quantity $\Lambda_n^\dagger(Q)$ has actually received significant
attention in its own right.  Thus 
Masser \cite{masser} has drawn upon the 
proof of \eqref{55} to show that
\beq\lab{98}
\Lambda_n^\dagger(Q) \ll \Mq^{n/2},
\eeq
for arbitrary indefinite quadratic forms
$Q$.   By adapting Kneser's example, Masser also shows that the exponent $n/2$ is best possible.  It is somewhat
surprising that the condition $x_1\neq 0$ inflates the corresponding exponent
of Cassels by $1/2$.  In this setting, our estimates for
$\Lambda_n^\dagger(Q)$ have rather more currency than those for 
$\Lambda_n(Q)$ did, since there is now no
analogue of the argument that was used to prove Theorem \ref{diag}.
In particular we can show that Masser's estimate is not best possible
for diagonal indefinite quadratic forms, as given by \eqref{Q}.
For such forms, the set of eigenvalues
$\{\lambda_1,\ldots,\lambda_n\}$ coincides precisely with the set of
coefficients $\{A_1,\ldots,A_n\}$.  We may therefore employ the
lower bound $n-3-\al_n\geq 2$, together with the fact that
$\mq^{-1}|\D| \leq \Mq^{n-1}$, in order to deduce the following
trivial consequence of Theorem \ref{c}.

\begin{cor}\lab{cor1}
Let $n\geq 5$ and assume that $Q \in \Z[X_1,\ldots,X_n]$ is a
diagonal indefinite quadratic form, with discriminant
$\D\neq 0$.  Then we have
$$
\Lambda_n^\dagger(Q) \ll_\ve  \Mq^{(2n(1+\be_Q)-1)/(n-3-\al_n)+\ve},
$$
for any $\ve>0$, where $\al_n$ is given by \eqref{parity} and $\be_Q$ is given by \eqref{8-q}.
\end{cor}

It is not hard to check that Corollary \ref{cor1} improves upon Masser's bound \eqref{98} when $n
\geq 7$ if $\D$ is odd and square-free, and when $n\geq 9$ in general.
In fact the improvement is quite substantial, suggesting a limiting exponent of $2$ as $n
$ gets large.

\subsection{Positive definite forms}

In this section we turn to the second major theme of this paper.  Given a
positive definite quadratic form $Q \in \Z[X_1,\ldots,X_n]$, we are now interested in determining whether or not
$\mcal{S}(k;Q)$ is empty, when $k$ is large.  An obvious necessary
condition for $\mcal{S}(k;Q)$ to be non-empty is that the congruence
$$
Q(\x)\equiv k \mod{p^t}
$$ 
should be soluble for every prime power $p^t$.  Let us say that the
pair $k,Q$ satisfies the ``weak local solubility condition''
if this occurs.  We shall usually just write  ``$(k,Q)$ satisfies 
weak \textsf{LSC}'',  for short. 
It is natural to question whether $\mcal{S}(k;Q)$ is
automatically non-empty for any pair $k,Q$ satisfying the weak 
local solubility condition.
The answer to this is negative, as demonstrated by the quadratic form
\beq \lab{eg1}
Q_1(\X)=2(X_1^2+\cdots+X_{n-1}^2)+(k+2)X_n^2.
\eeq
This example is due to Watson \cite[\S 7.7]{W'}. It is not hard to check that 
$(k,Q_1)$ satisfies weak \textsf{LSC}, but that the equation $Q_1(\x)=k$ is
insoluble for large odd values of $k$.  
It is at this point that the work of  Tartakowsky \cite{tart} enters the picture.
Define
$$
\K_n(Q):=\{k\in\N:~ \mcal{S}(k;Q)=\emptyset, ~\mbox{$(k,Q)$ satisfies
  weak \textsf{LSC}}\},
$$
and 
$$
\k_n(Q):=\max_{k\in\K_n(Q)} \{k\}.
$$
Then under the assumption that $n \geq 5$, he has shown that $\K_n(Q)$
is finite.  

Tartakowsky's argument does not lead to any estimate for
the cardinality of $\K_n(Q)$, and the problem of finding an effective upper
bound for this quantity has since been considered by several authors. 
One of the most impressive results in this direction is due to Watson
\cite{W}, giving 
that 
\beq\lab{8-watson}
\k_n(Q)\ll
\left\{
\begin{array}{ll}
|\D|^{5/(n-4)+1/n}, & \mbox{if $5 \leq n \leq 9$,}\\
|\D|, & \mbox{if $n \geq 10$.}
\end{array}
\right.
\eeq
Note that a positive definite quadratic form is automatically
non-singular, and furthermore, we trivially have $\#\K_n(Q)\leq \k_n(Q)$.
It should be clear from \eqref{eg1} that Watson's bound is 
best possible for $n\geq 10$, as we have $k \gg \Delta_{Q_1}$ here.
Hsia and Icaza \cite[\S 4]{h-i} have since provided the estimate
\beq\lab{8-hi}
\k_n(Q)\ll 
|\D|^{(n-2)/(n-4)+2/n},
\eeq
for $n\geq 5$, in which the implied constant is made completely explicit.
This is sharper than \eqref{8-watson} when $n=5$ or $n=6$.

In the intermediate case $5 \leq n \leq 9$, Watson obtains
sharper bounds in the special case of diagonal quadratic forms. In fact his approach to the problem is to
first handle the case of diagonal forms through a classical
application of the circle method, before then combining these
results with a diagonalisation process to handle the general case.  
This is somewhat wasteful and our approach to the problem will involve
handling the case of general quadratic forms directly.
It transpires that our method is most effective when the height
$\Mq$ of $Q$ is small compared to $\D$, whereas Watson's method is best
when $\Mq$ is large compared to $\D$.  By merging the two approaches we shall
succeed in beating Watson's bound when $5 \leq n\leq 9$.  

\begin{theorem}\label{iterate}
Let $n\geq 5$ and assume that $Q \in \Z[X_1,\ldots,X_n]$ is a
positive definite quadratic form.  Then we have
$$
\k_n(Q)\ll_\ve
\left\{
\begin{array}{ll}
|\D|^{\phi(n)+\ve}, & \mbox{if $5 \leq n \leq 9$,}\\
|\D|, & \mbox{if $n \geq 10$,}
\end{array}
\right.
$$
for any $\ve>0$, where
\beq\lab{phi}
\phi(n)=\frac{4(n-2)(3n^2-7n-3)}{(2n^3-9n^2+2n+12)(n-3)}.
\eeq
\end{theorem}

In order to facilitate comparison between Theorem \ref{iterate} and
the bounds in \eqref{8-watson} and \eqref{8-hi}, we have calculated approximate values for the
exponents in the following table:

\begin{center}
\begin{tabular}{|l|| l| l|l|}
\hline
$n$ & $\frac{5}{n-4}+\frac{1}{n}$ & $\frac{n-2}{n-4}+\frac{2}{n}$ & $\phi(n)$\\
\hline
\hline
$5$ & $5.200..$ & $3.400..$ & $4.723..$\\
$6$ & $2.666..$ & $2.333..$ & $2.545..$\\
$7$ & $1.809..$ & $1.952..$ & $1.752..$\\
$8$ & $1.375..$ & $1.750..$ & $1.341..$\\
$9$ & $1.111..$ & $1.622..$ & $1.088..$\\
\hline
\end{tabular}
\end{center}

One easily checks 
that both of the exponents $5/(n-4)+1/n$ and $\phi(n)$ are
strictly less than $1$ for $n \geq 10$. Moreover the bound
of Hsia and Icaza is the best available for $n=5$ and $n=6$, but is 
weaker than Theorem \ref{iterate} for $7\leq n \leq 9$.
Inspired by our results in the previous section, it might be expected
that sharper bounds are available for $\k_n(Q)$ when the height 
of $Q$ is not too large compared to the discriminant. The following
result shows that this is indeed the case.

\begin{theorem}\label{8-w'}
Let $n\geq 5$ and assume that $Q \in \Z[X_1,\ldots,X_n]$ is a
positive definite quadratic form.  Then we have
$$
\k_n(Q)\ll_{\ve} \big(|\D|^{(n-2)/(n-4)}\Mq^{n+\ve}\big)^{2/(n-3)},
$$
for any $\ve>0$.
\end{theorem}

Let $\Sigma_1\leq \cdots \leq \Sigma_n$ denote the $n$ successive
minima of $Q$. If $Q$ is Minkowski reduced, in the sense of 
Watson \cite[\S 2.9]{W'}, for example, then it follows that
$$
\Mq \ll \Sigma_n.
$$
In view of the fact that the set $\K_n(Q)$ is left
invariant under any unimodular transformation,
it is not hard to see that the statement of Theorem
\ref{8-w'} remains true with $\Mq$ replaced by $\Sigma_n$.
This gives a version of the result that is independent of the
particular choice of coordinates.
For a typical quadratic form $Q$ one expects $\Sigma_n$ to have order of magnitude $|\D|^{1/n}$.
This yields
$$
\k_n(Q)\ll_{\ve} |\D|^{4/(n-4)+\ve}
$$ 
in  Theorem \ref{8-w'}. Not only does this improve on \eqref{8-watson}
for every $n\geq 5$, it does so by quite an ample margin.

It is natural to ask about the corresponding situation for smaller
values of $n$.  When $n=3$, Duke \cite{duke} has used the theory of weight
$\frac{3}{2}$ modular forms to tackle the problem.  Specifically, he
has shown that there exists an absolute constant $c>0$ such that
if $k$ is square-free, with $k>c|\D|^{337}$, then 
the equation $Q(X_1,X_2,X_3)=k$ has an
integral solution provided that $(k,Q)$ satisfies weak \textsf{LSC}.
Here, the constant $c$ is ineffective, since it arises out of an application of Siegel's lower bound for $L(1,\chi)$.

When $n=4$, it turns out that stronger assumptions are
needed to ensure the solubility of $Q(X_1,X_2,X_3,X_4)=k$.
Consider the positive definite quadratic form
$$
Q_2(\X)=X_1^2+X_2^2+7(X_3^2+X_4^2).
$$
Watson \cite[\S 7.7]{W'} has observed that by taking $k=3\times
7^{2u}$, for $u \in \N$, one can show that the set $\K_4(Q_2)$ is
infinite.  Hence a stronger local solubility condition is needed to treat the case $n=4$. Let 
\beq\lab{sangram}
\tau_p=\left\{
\begin{array}{ll}
0, & p>2,\\
1, & p=2,
\end{array}
\right.
\eeq
and let $Q\in\Z[X_1,\ldots,X_n]$ be a positive definite
quadratic form.  Then we shall say that the pair $k, Q$ 
satisfies the ``strong local solubility condition'', usually written ``$(k,Q)$
satisfies strong \textsf{LSC}'' for short, if for every prime $p$
there exists $\x\in(\Z/p^{1+2\tau_p}\Z)^n$ such that
\begin{equation}\lab{ganga}
Q(\x)\equiv k \mod{p^{1+2\tau_p}}, \quad p\nmid \ma{A}\x.
\end{equation}
We shall occasionally say that ``$(k,Q)$
satisfies strong \textsf{LSC} modulo $p$'' if this congruence has a
solution for that particular choice of prime $p$.
We now define the set
$$
\K_n^*(Q):=\{k\in\N:~ \mcal{S}(k;Q)=\emptyset, ~\mbox{$(k,Q)$ satisfies
  strong \textsf{LSC}}\},
$$
and the corresponding quantity
$$
\k_n^*(Q):=\max_{k\in\K_n^*(Q)} \{k\}.
$$
Clearly $\k_n^*(Q)$ is finite for $n\geq 5$, since $\K_n^*(Q)\subseteq
\K_n(Q)$.  The corresponding result for $n=4$ is 
due to Fomenko \cite{fom} and uses the theory of modular forms. 
Schulze-Pillot \cite{s-p} has since refined the argument, obtaining
$$
\k_4'(Q) \ll_{\ve} N^{14+\ve},
$$
for a quantity $\k_4'(Q)$ similar to $\k_4^*(Q)$,
where $N$ denotes the level of $Q$. 
%T%REMOVE THE DEFN OF LEVEL AND OUR EARLIER REMARKS..
%This is defined to be the least positive
%integer $N$ for which $N\ma{H}^{-1}$ has integer coefficients, where
%$\ma{H}$ denotes the Hessian matrix associated to $Q$. 
%The level of $Q$ is intimately related to the discriminant
%of $Q$, although the precise connection is still not well understood.
It should be noted that Schulze-Pillot's bound is completely explicit,
and that he achieves finer estimates under the assumption that $N$ is square-free.
The essential difference between $\k_4'(Q)$ and $\k_4^*(Q)$ is that
there should be primitive local solutions everywhere. This 
is implied by our strong local solubility conditions 
\eqref{ganga}, whence $\k_4^*(Q)\leq \k_4 '(Q)$. 
Hanke \cite[Theorem 6.3]{hanke} has also used a modular forms interpretation to
examine a quantity similar to $\k_4^*(Q)$,
but the estimate he arrives at is too complicated to state
here. Again, an alternative local solubility condition is 
employed, which differs from both Schulze-Pillot's and ours. It
corresponds to assuming weak \textsf{LSC}, together with an extra local condition on $k$
for those primes $p$ such that $Q$ is anisotropic modulo $p$. 
Our approach leads us to the following result.

\begin{theorem}\label{w}
Let $n\geq 4$ and assume that $Q \in \Z[X_1,\ldots,X_n]$ is a
positive definite quadratic form.  Then we have
$$
\k_n^*(Q)\ll_{\ve} \big(|\D|\Mq^{n+\ve}\big)^{2/(n-3)},
$$
for any $\ve>0$. 
\end{theorem}

Take $n=4$ in the statement of Theorem \ref{w}.
Then, using the fact that the set $\K_n^*(Q)$ is left
invariant under any unimodular transformation,
together with the basic property $\Mq\ll |\D|$ satisfied by any 
Minkowski reduced quadratic form $Q$,
we deduce the following result.

\begin{cor}
Let $Q \in \Z[X_1,X_2,X_3,X_4]$ be a
positive definite quadratic form.  Then we have
$$
\k_4^*(Q) \ll_{\ve}|\D|^{10+\ve},
$$
for any $\ve>0$. 
\end{cor}

Returning to the generic setting, for which $\Mq^n \ll |\D|
\ll \Mq^{n}$, it follows from Theorem \ref{w} that 
$\k_n^*(Q)\ll_{\ve} |\D|^{4/(n-3)+\ve}$
for typical positive definite quadratic forms in $n\geq 4$
variables. In particular we have the sharper bound 
$\k_4^*(Q)\ll_{\ve} |\D|^{4+\ve}$ for generic quaternary forms $Q$.

\subsection{Outline of the paper}

The underlying tool in this paper is a modern form of the  Hardy--Littlewood circle
method, due to Heath-Brown \cite{HB'}. This will be discussed
in more detail in the following section. 
In \S \ref{pf-iterate}, which is essentially independent of the circle method, we
shall undertake the proof of Theorem \ref{iterate}.
This part of the paper involves a delicate reduction argument which allows one to study the
equation $Q=k$ under stronger local solubility assumptions, thereby
permitting an application of Theorem~\ref{w}, rather than the obvious
application of Theorem \ref{8-w'}.
Once combined with Watson's proof of \eqref{8-watson}, this will be
enough to furnish the statement of Theorem \ref{iterate}.

Our use of the circle method begins in earnest in \S \ref{act}. 
In fact we shall use it to establish
an asymptotic formula for the number of $\x \in \Z^n$ such that
$Q(\x)=k,$ which are constrained to lie in a certain expanding region.
A crucial feature of our asymptotic formula is that its dependence
upon the coefficients of $Q$ needs to be made 
completely explicit. 
Indeed, we will then be able to determine precise 
conditions on the size of the region (resp. the size of $k$) 
needed to ensure that $\mcal{S}(k;Q)$ is non-empty, as required for
Theorem \ref{c} (resp. Theorems \ref{8-w'} and \ref{w}). 
This requires an appreciable amount of work,
since the formulation of the method given by Heath-Brown  \cite{HB'} pays no
attention to the question of uniformity in the coefficients of the
quadratic form.   

A further obstacle that we'll need to deal with, and which marks
another departure from the usual applications of the circle method, emerges in the 
treatment of the singular series $\ss(k,Q)$. Thus we will need to
bound $\ss(k,Q)$ away from zero uniformly in terms of $k$ and the coefficients of
$Q$. While the obvious approach for doing this
would undoubtedly give something here, we have adopted a rather more
sophisticated argument in \S \ref{ss}, and the bounds obtained are
actually quite sharp.

\begin{notation}
Throughout this paper we shall write $\int f(\z) \d\z$ for
the $n$-fold repeated integral of $f(\z)$ over $\R^n$.  Given $q \in
\N$, a sum with a condition of the form $\ma{b}\tmod q$ will mean
a sum taken over $\ma{b} \in \Z^n$ such that the components of
$\ma{b}$ run from $0$ to $q-1$. Finally, for any $\al \in \R$ we shall
write $e(\al):=e^{2\pi i \al}$ and $e_q(\al):=e^{2\pi i \al/q}$.
\end{notation}

\begin{ack}
The authors are grateful to Roger Heath-Brown for several
useful conversations relating to the contents of this
paper. 
\end{ack}

\section{Preliminaries}\lab{2}

In this section we shall collect together the main ingredients in the proofs
of Theorems~\ref{c}, \ref{8-w'} and \ref{w}. Let $n \geq 4$ and let $Q\in \Z[X_1,\ldots,X_n]$ 
be a non-singular quadratic form of discriminant $\D$, with
underlying matrix $\mathbf{A}$. 
Let $\mathbf{R} \in \mathrm{SO}_n(\R)$ be an orthogonal matrix
that diagonalises $\mathbf{A}$, which we regard as being fixed throughout this paper.
Then there exist $\lambda_1,\ldots,\lambda_n \in \R$,
such that 
\beq\lab{choice-R}
\mathbf{R}^T\mathbf{A} \mathbf{R}=\mathrm{Diag}(\lambda_1,\ldots,\lambda_n),
\eeq
and $\lambda_1\cdots \lambda_n =\D$.  
In particular, we have $Q(\mathbf{R}\ma{U})=\lambda_1U_1^2+\cdots +\lambda_n
U_n^2$, and since $\lambda_1,\ldots,\lambda_n$ are merely the eigenvalues
of $\mathbf{A}$, \eqref{eigen} clearly holds.

Given an arbitrary polynomial $f\in \Z[X_1,\ldots,X_n]$, 
and a bounded function $w: \R^n \rightarrow \R_{\geq 0}$
of compact support, we define the weighted counting function
$$
N_w(f;B):=\sum w(B^{-1}\x), 
$$
for any $B \geq 1$. Here, the summation is taken over all $\x \in \Z^n$ for which $f(\x)=0$.
Our proof of Theorem \ref{c} is based upon an analysis of the asymptotic
behaviour of $N_w(Q;B)$, as $B \rightarrow \infty$, for a suitable $w$.  
Likewise, to prove Theorems \ref{8-w'} and \ref{w} we shall study the
counting function $N_w(Q-k;k^{1/2})$, as $k \rightarrow \infty$.
The quantities $N_w(Q;B)$ and $N_w(Q-k;k^{1/2})$ have received considerable attention over the years,
and several methods have been developed to study them.  The method that we
shall employ is based upon the new form of the Hardy--Littlewood circle
method, due to Heath-Brown \cite{HB'}.  
Using this version of the circle method, Heath-Brown has established the existence of a
non-negative constant $c_w(k,Q)$ such that 
$$
N_w(Q-k;k^{1/2})=c_w(k,Q)k^{n/2-1} + O_{\ve,Q}\big(k^{(n-1)/4+\ve}\big),
$$
provided that $n\geq 4$ and $w$ belongs to a certain class of
weight functions.  Similarly, for $n\geq 5$ and the same class of
weights $w$, he shows that there is a non-negative constant $c_w(Q)$ such that 
$$
N_w(Q;B)=c_w(Q)B^{n-2}+ 
\left\{
\begin{array}{ll}
O_{\ve,Q}(B^{n/2+\ve}), & \mbox{if $n\geq 6$ is even,}\\
O_{\ve,Q}(B^{(n-1)/2+\ve}), & \mbox{if $n\geq 5$ is odd.}
\end{array}
\right.
$$ 
In order to discuss the two cases simultaneously, it will be
convenient to think of the case $k=0$ as corresponding to a study of 
$N_w(Q;B)$, as $B\rightarrow \infty$, and the case $k>0$ as
corresponding to a study of $N_w(Q-k;k^{1/2})$, as $k \rightarrow
\infty$. In either case we may therefore refer to the counting
function $N_w(Q-k;B)$, for $k \geq 0$, with the understanding that
$B=k^{1/2}$ when we are in the case $k>0$. On assuming that $n\geq 5$
in the case $k=0$, we may combine Heath-Brown's estimates to deduce
that 
\beq\lab{hb}
N_w(Q-k;B)=c_w(k,Q)B^{n-2} + O_{\ve,Q}\big(B^{(n-1+\gamma_n)/2+\ve}\big),
\eeq
for a suitable constant $c_w(k,Q)\geq 0$, where
\beq\lab{even-odd}
\gamma_n=
\left\{\begin{array}{ll} 1, & \mbox{if $n$ is even and $k=0$,}\\ 
0, & \mbox{otherwise.}
\end{array}
\right.
\eeq
As indicated above, the central component in our work is a finer version of \eqref{hb},
in which the dependence upon the coefficients of $Q$ is made completely explicit.
On establishing a suitable lower bound for the constant $c_w(k,Q)$, it will
then be possible to determine precise information about the size of
$B$ that is needed to ensure that $N_w(Q-k;B)>0$.

The constant $c_w(k,Q)$ may be interpreted as a product of local
densities, and we proceed to discuss it in more detail.
For any prime $p$, the $p$-adic density of solutions is defined to be
\beq\lab{sig-p}
\sigma_p=\sigma_p(k,Q):=\lim_{t \rightarrow \infty} p^{-t(n-1)}
N(p^t),
\eeq
where
\beq\lab{N}
N(p^t):=\#\{\x \mod{p^t}: Q(\x)\equiv k \mod{p^t}\}.
\eeq
When this limit exists, the singular series is given by 
\beq\lab{singseries}
\ss(k,Q):=\prod_{p}\sigma_p.
\eeq
We shall set $\ss(Q):=\ss(0,Q)$.
It transpires that $\ss(k,Q)$ is always convergent for the quadratic forms considered here.
Let us write
\beq\lab{N*}
N^*(p^t):=\#\{\x \mod{p^t}: ~Q(\x)\equiv k \mod{p^t}, ~p \nmid \ma{A}\x\},
\eeq
for any prime power $p^t$.  It follows from a simple application of Hensel's lemma that
\beq\lab{hensel}
N^*(p^t)\geq p^{(n-1)(t-1-2\tau_p)}N^*(p^{1+2\tau_p}),
\eeq
for any $t \geq  1+2\tau_p$, where $\tau_p$ is given by \eqref{sangram}.
We shall make use of this inequality at several points of our argument.

In order to introduce the singular integral, it will be
convenient to define the polynomial $P_Q\in \Z[X_1,\ldots,X_n]$ according
to the rule
\beq\lab{P}
P_Q(\X):=\left\{
\begin{array}{ll}
Q(\X), & \mbox{if $k=0$,}\\ 
Q(\X)-1, & \mbox{if $k>0$.} 
\end{array}
\right.
\eeq
In particular we clearly have $P_Q(\X)=B^{-2}(Q(B\X)-k)$, whether or not
$k$ is zero.
Then for any infinitely differentiable function $w: \R^n
\rightarrow \R_{\geq 0}$ of compact support $\supp(w)$, such that $\nabla P_Q\neq \ma{0}$ on the closure
of $\supp(w)$, the  corresponding singular integral is defined to be 
\beq\lab{singint-1}
\sigma_\infty(w;P_Q):=\lim_{\ve\rightarrow 0}(2\ve)^{-1}\int_{|P_Q(\x)|\leq
\ve} w(\x)\d\x.
\eeq
This limit exists, and moreover is positive if $w$ takes a positive value for
some real solution $\x$ of $P_Q(\x)=0$, by the first part of
\cite[Theorem 3]{HB'}.
One should think of $\sigma_\infty(w;P_Q)$ as giving the real
density of solutions, weighted by $w$.  With these definitions in
mind, we then have  $c_w(k,Q)=\sigma_\infty(w;P_Q)\ss(k,Q)$ in
\eqref{hb}.

Before revealing our uniform version of \eqref{hb}, we must first
decide upon the choice of weight function that we shall work with.
Consider the function $w_0: \R \rightarrow \R_{\geq 0}$, given by
\beq\lab{w0} 
w_0(x) := \left\{ \begin{array}{ll}  e^{-(1-x^2)^{-1}}, &
\mbox{if $|x| < 1$},\\
0, & \mbox{if $|x| \ge 1$}.
\end{array} 
\right.
\eeq 
Then $w_0$ is infinitely differentiable with compact support
$[-1,1]$.  Now define the function
\beq\lab{w1}
w_1(\x):=w_0(2x_1-2)w_0(x_2)\cdots w_0(x_n),
\eeq
on $\R^n$.  Then $w_1$ is infinitely differentiable, with
support $[\frac{1}{2},\frac{3}{2}]\times [-1,1]^{n-1}$.
Recall the orthogonal matrix $\mathbf{R}\in \mathrm{SO}_n(\R)$ that was chosen so that
\eqref{choice-R} holds. Then 
we shall work with the function $w_Q:\R^n\rightarrow \R_{\geq 0}$,
given by 
\beq\lab{wQ}
w_Q(\x):=\tilde{w}(\mathbf{R}^T \x),
\eeq
where 
\beq\lab{tW}
\tilde{w}(\x):=w_1(|\lambda_1|^{1/2}x_1, \ldots, |\lambda_n|^{1/2}x_n).
\eeq
Let us write 
\beq\lab{R-sig}
Q_{\mathrm{sgn}}(\X):=\sigma_1 X_1^2+\sigma_2 X_2^2+\cdots+\sigma_nX_n^2,
\eeq
where $\sigma_i:=\lambda_i/|\lambda_{i}|$ for $1\leq i \leq n$.    Then 
 \eqref{singint-1} implies that 
\begin{align*}
\sigma_\infty(w_Q;P_Q)
&=
\lim_{\ve\rightarrow 0}(2\ve)^{-1}\int_{|P_Q(\x)|\leq
\ve} \tilde{w}(\mathbf{R}^T \x)\d\x\\
&=
\lim_{\ve\rightarrow 0}(2\ve)^{-1}\int_{|P_Q(\mathbf{R}\u)|\leq
\ve} w_1(|\lambda_1|^{1/2}u_1, \ldots, |\lambda_n|^{1/2}u_n)\d\u\\
&=\frac{1}{|\D|^{1/2}}\lim_{\ve\rightarrow 0}(2\ve)^{-1}\int_{|P_{Q_{\mathrm{sgn}}}(\ma{v})|\leq
\ve} w_1(\ma{v})\d\ma{v} =\frac{\sigma_\infty(w_1;P_{Q_{\mathrm{sgn}}})}{|\D|^{1/2}}.
\end{align*}
We shall write 
\beq\lab{sig-inf}
\sigma_{\infty}=\sigma_{\infty}(w_1;P_{Q_{\mathrm{sgn}}}),
\eeq
for convenience, where $Q_{\mathrm{sgn}}$ is given by \eqref{R-sig} and 
$P_{Q_{\mathrm{sgn}}}$ is given by \eqref{P}. In particular it follows from our remarks
above that 
\beq\lab{mendip3}
1\ll \sigma_\infty \ll 1.
\eeq
We are now ready to reveal the main ingredient in our work. The
following result will be established in \S \ref{act}.

\begin{pro}\lab{main'}
Let $n\geq 4$ and $k>0$, or $n\geq 5$ and $k=0$. Then we have
$$
N_{w_Q}(Q-k;B) = \frac{\sigma_\infty \ss(k,Q) }{|\D|^{1/2}} B^{n-2}+
O_{\ve}\big( 
\Mq^{n/2+\ve}B^{(n-1+\gamma_n)/2+\ve}\big),
$$
where $\gamma_n$ is given by \eqref{even-odd}.
\end{pro}

As a method for proving results of the sort in Theorems \ref{c}, \ref{8-w'} and
\ref{w}, Proposition~\ref{main'} is not altogether new. 
In fact the second author \cite[Theorem 2]{dietmann} has established an asymptotic formula for
a quantity very similar to $N_{w_Q}(Q-k;B)$. He does so for precisely the same ranges of
$n$, and also obtains uniformity with respect to the coefficients of $Q$.
However the error term that we obtain is
substantially sharper than that obtained there.
%Although we shall not pursue it here, it is possible to improve upon
%the error term in Proposition \ref{main'} when $Q$ is diagonal.
%This would allow us to replace the
%term $\Mq^{n/2}$ by $\Mq^{(n-1)/2}$, which would ultimately lead to
%a version Corollary \ref{cor1} with exponent
%$(2n-2)/(n-3-\al_n)$, instead of $(2n-1)/(n-3-\al_n)$.
%R% better in a further paper; application to vector sieve for n=4!
%T%OK!
It is worth highlighting that the classical form of the
circle method (see Davenport \cite{dav}, for
example) could also be used to establish a result of
the type in Proposition~\ref{main'} when $n\geq 5$. However, a 
single Kloosterman refinement is needed to treat the case $n=4$ and
$k>0$.

The weight $w_Q$ that occurs in Proposition \ref{main'} has been
specially chosen to optimise the error
term in the asymptotic formula. When $k=0$ it is of independent interest to
try and obtain versions of this result for the counting function
$N_w(Q;B)$ associated to a weight $w:\R^n\rightarrow\R_{\geq 0}$ that
approximates the characteristic function of  $[-1,1]^n$, since
this amounts to counting rational points of bounded height on the
quadric hypersurface $Q=0$ in $\bfP^{n-1}$.
This line of enquiry has been pursued by the first author
\cite{qupper} for diagonal quadratic forms. A novel feature of this
work is that quaternary forms are handled, these not being 
touched upon in the present work when $k=0$.  

Returning to Proposition \ref{main'}, it 
is clear that we shall also need
some control over the size of the singular series $\ss(k,Q)$ appearing
in Proposition \ref{main'}.  We shall be able to do so 
under the assumption that $(k,Q)$ satisfies weak \textsf{LSC}, and we
shall get our sharpest bound when it is assumed that 
$(k,Q)$ satisfies strong \textsf{LSC}.
The following result will be established in \S \ref{ss}.

\begin{pro}\lab{ss'}
Let $\ve>0$.  Suppose that $n\geq 5$ and 
$(k,Q)$ satisfies weak \textsf{LSC}, with $Q$ a 
non-singular quadratic form. Then 
$$
\ss(k,Q) \gg_{\ve} |\D|^{-\theta_{k,Q}-\ve},
$$
where
\beq\lab{8-theta}
\theta_{k,Q}=
\left\{
\begin{array}{ll}
0, & \mbox{if $(k,Q)$ satisfies strong \textsf{LSC},}\\
1/(n-4), &\mbox{otherwise}.
\end{array}  
\right.
\eeq
Now suppose that $n=4$ and $(k,Q)$ satisfies strong \textsf{LSC}, with $k>0$.
Then 
$$
\ss(k,Q) \gg_{\ve} k^{-\ve}|\D|^{-\ve}.
$$
\end{pro}

It is now an easy matter to combine Propositions \ref{main'} and
\ref{ss'} to deduce Theorems~\ref{c},  \ref{8-w'} and \ref{w}.   
The deduction of Theorem \ref{iterate} is rather more involved and will
be undertaken in the subsequent section.   
Let us begin by deriving Theorem \ref{c}, for which we shall apply Propositions 
\ref{main'} and  \ref{ss'} in the case $n\geq 5$ and $k=0$.
In particular we have $\gamma_n=\al_n$ in Proposition \ref{main'}, where
$\al_n$ is given by \eqref{parity}.
When the discriminant $\D$ of $Q$ is square-free, an application of
the Chevalley--Warning theorem implies that $(0,Q)$ satisfies strong \textsf{LSC} modulo $p$, for every
odd prime $p$.
%R%Yes, but for p=2 we have a problem: We need there, say in the diagonal
%R%case, at least five odd coefficients. For example, x_1^2+...+x_4^2
%R%does not satisfy strong LSC for k=0! But for n \ge 5 it is sufficient
%R%that the determinant is odd.
If $n \ge 5$ and $\D$ is odd, then it is easily seen that
$(0,Q)$ satisfies strong \textsf{LSC} modulo $2$.
Thus we may take $\theta_{0,Q}=\beta_Q$ in our application of Proposition
\ref{ss'}, where $\beta_Q$ is given by \eqref{8-q}.
We may therefore combine Propositions~\ref{main'} and  \ref{ss'} with 
\eqref{mendip3}, in order to deduce that $N_{w_Q}(Q;B)>0$ provided
that $B$ is chosen so that
$$
B\gg_{\ve} \big(|\D|^{1+2\beta_Q} \Mq^{n+\ve}\big)^{1/(n-3-\al_n)}.
$$
Recall that $\lambda_1,\ldots,\lambda_n\in \R$ denote the eigenvalues of the underlying
matrix $\ma{A}$.  Then it follows that there is at least one non-trivial vector $\x \in \Z^n$,
for which $Q(\x)=0$ and $x_1\neq 0$, with $|\x|\leq \tilde{B}$, provided that 
$$
\tilde{B}
\gg_{\ve} \max_{1\leq i\leq n}\{|\lambda_i|^{-1/2}\}
\big(|\D|^{1+2\be_Q} \Mq^{n+\ve}\big)^{1/(n-3-\al_n)}.
$$
On recalling the definition \eqref{min-Q} of
$\mq$, this therefore completes the proof of Theorem \ref{c}.

We conclude this section by deducing the statements of Theorems
\ref{8-w'} and \ref{w}.
For this we shall apply Propositions 
\ref{main'} and  \ref{ss'} in the case $n\geq 4$ and $k>0$.
In particular we have $B=k^{1/2}$ and $\gamma_n=0$ in Proposition \ref{main'}.
Let $Q\in \Z[X_1,\ldots,X_n]$ be a
positive definite quadratic form and suppose that $(k,Q)$ satisfies
weak \textsf{LSC}.  Then it follows from 
Propositions  \ref{main'} and \ref{ss'}, together with
\eqref{mendip3}, that $N_{w_Q}(Q-k;k^{1/2})>0$ provided
that 
$$
k^{1-\ve} \gg_{\ve} \big(|\D|^{1+2\theta_{k,Q}} \Mq^{n+\ve}\big)^{2/(n-3)},
$$
and $(k,Q)$ is assumed to satisfy strong \textsf{LSC} in the case $n=4$.
Here, $\theta_{k,Q}=0$ if $(k,Q)$ satisfies
strong \textsf{LSC}, and $\theta_{k,Q}=1/(n-4)$ otherwise.
This provides the required upper bounds for $\k_n(Q)$ and $\k_n^*(Q)$.

\section{A hybrid approach to Theorem \ref{iterate}}\lab{pf-iterate}

The purpose of this section is to deduce the statement of Theorem
\ref{iterate} from Theorem \ref{w}.
Throughout this section let $n\geq 5$, and let $Q \in \Z[X_1,\ldots,X_n]$
be a positive definite quadratic form of discriminant $\D$.  As we
have already mentioned, the sets $\K_n(Q)$ and $\K_n^*(Q)$ are left 
invariant under any unimodular transformation, since
$\mcal{S}(k;Q)$ is non-empty if and only if $\mcal{S}(k;Q')$ is
non-empty, for any $Q'$ that is equivalent to $Q$.  

We begin by recording a simple calculation for the number of solutions
to a quadratic congruence modulo an odd prime $p$.  Given integers
$k,a_1,\ldots,a_r$, and any odd prime $p$, we define
\beq\lab{mr-def}
M_r(p):=\#\{\z \mod{p}: ~a_1z_1^2+\cdots+a_rz_r^2\equiv k \mod{p}, ~p\nmid \z\},
\eeq
where $\z=(z_1,\ldots,z_r).$
The following result ought to be well-known, but we have included our
own proof for the sake of completeness.

\begin{lemma}\lab{mr}
Suppose that $p\nmid 2a_1\cdots a_r$. Then we have
$$
M_r(p)=\left\{
\begin{array}{ll}
p^{r-1}-\kappa_p + \Big(\frac{a_1\cdots a_r}{p}\Big)\omega_p^r p^{-1}(\kappa_p
p -1), & \mbox{if $r$ is even},\\
p^{r-1}-\kappa_p + \Big(\frac{-ka_1\cdots a_r}{p}\Big)\omega_p^{r+1} p^{-1}, & \mbox{if $r$ is odd},
\end{array}
\right.
$$
where $\omega_p:=i^{(p-1)^2/4}\sqrt{p}$ and 
\beq\lab{28-kappa}
\kappa_p:=\left\{
\begin{array}{ll}
1, & \mbox{if $p \mid k$}, \\
0, & \mbox{if $p\nmid k$}.
\end{array}
\right.
\eeq
\end{lemma}

\begin{proof}
In order to study $M_r(p)$, we write
\begin{align*}
M_r(p)&=-\kappa_p+{p^{-1}}\sum_{c\mod p} \sum_{\z\mod{p}} e_p\big(c(a_1z_1^2+\cdots +a_rz_r^2-k)\big)\\
&= p^{r-1}-\kappa_p + p^{-1} \sum_{c=1}^{p-1} e_p(-ck)\prod_{i=1}^r
\sum_{z \mod{p}} e_p(c a_i z^2).
\end{align*}
Since $p \nmid 2ca_1\cdots a_r$, the innermost sum is a Gauss sum and
so takes the value $\big(\frac{ca_i}{p}\big)\omega_p$.  Thus we obtain
$$
M_r(p)=p^{r-1}-\kappa_p + \Big(\frac{a_1\cdots a_r}{p}\Big)\omega_p^r
p^{-1}\sum_{c=1}^{p-1} \Big(\frac{c}{p}\Big)^r e_p(-kc). 
$$
Once combined with the well-known equalities
$$
\sum_{c=1}^{p-1} \Big(\frac{c}{p}\Big)^r e_p(ac)=\left\{
\begin{array}{ll}
p-1, & \mbox{if $r$ is even and $p\mid a$,}\\
-1, & \mbox{if $r$ is even and $p\nmid a$,}\\
(\frac{a}{p})\omega_p, & \mbox{if $r$ is odd,}
\end{array}
\right.
$$
we easily conclude the proof of Lemma \ref{mr}.
\end{proof}

We are now ready to commence the proof of Theorem \ref{iterate}.
Our first step is a certain reduction argument that will render it sufficient to
examine the solubility of the equation $Q=k$ under stronger local
solubility assumptions.  Recall the statement \eqref{ganga} of strong
\textsf{LSC} modulo a prime $p$, and the definition \eqref{sig-p} of
$\sigma_p(k,Q)$. Then we have the following result.

\begin{lemma}\lab{reduction}
Assume that $n \geq 5$ and $(k, Q)$ satisfies weak \textsf{LSC}.
Then there exists a positive integer $k'\leq k$ and a
positive definite quadratic form
$Q'\in\Z[X_1,\ldots,X_n]$, such that the following hold:
\begin{enumerate}
\item[(i)]
$(k',Q')$ satisfies  strong \textsf{LSC} modulo every $p>2$, and $\sigma_2(k',Q') \gg 1$.
\item[(ii)]
$\mcal{S}(k;Q)= \emptyset$ if and only if $\mcal{S}(k';Q')= \emptyset$.
\item[(iii)]
$k'/|\Delta_{Q'}| \ge k/|\Delta_Q|$ and $|\Delta_{Q'}|\leq |\D|$.
\end{enumerate}
\end{lemma}

\begin{proof}
If $(k, Q)$ satisfies strong \textsf{LSC}, then we
may set $Q'=Q$ and $k'=k$.  Indeed, it follows from \eqref{sig-p} and
\eqref{hensel} that 
$$
\sigma_2(k,Q)\geq 2^{-3(n-1)}\gg 1,
$$
if $(k, Q)$ satisfies strong
\textsf{LSC} modulo $2$.  
Otherwise, if $(k,Q)$ fails to satisfy strong \textsf{LSC}, then our
goal will be to derive the existence of a positive integer $k'<k$,
and a positive definite quadratic form
$
Q'\in\Z[X_1,\ldots,X_n],
$ 
such that $(k',Q')$ satisfies weak \textsf{LSC}, and conditions (ii),
(iii) are satisfied.  Once this is achieved it is clear how to
complete the proof of Lemma \ref{reduction}: either 
we can show that the pair $k',Q'$ satisfies condition (i), in which case we are done,
or else we may iterate the argument to produce a new pair $k'',Q''$.
That this process must terminate after a finite number of steps is obvious from the fact that 
$1 \leq \cdots < k'' <k'<k$.   

Let us begin by supposing that the pair $(k,Q)$ does not satisfy strong
\textsf{LSC} modulo $p$, for an odd prime $p$.  In this setting
it is well known that $Q$ can be diagonalised over the ring $\Z_p$ of
$p$-adic integers (see \cite[\S 4.3]{W'}, for example).  Since
this process does not alter whether or not the pair $k,Q$ satisfies strong \textsf{LSC}, we may assume
that $Q(\X) =A_1 X_1^2 + \cdots + A_n X_n^2$. After a  change
of variables we may further assume that 
\beq\lab{change}
Q(\X)\equiv a_1X_1^2+\cdots+a_rX_r^2 \mod{p}, 
\eeq
for some $0\leq r \leq n$, with $p \nmid a_1\cdots a_r$. 
If $r\geq 3$ then Lemma~\ref{mr} implies
that the congruence  $Q(\x)\equiv k \tmod{p}$ has a solution with
$p \nmid A_ix_i$ for some $1\leq i \leq n$.  
This gives a solution of \eqref{ganga}, which is contrary to our assumption. 
Thus we may assume that $r \leq 2$ in \eqref{change}.  
We claim that the only possibility is  $\kappa_p=1$, where $\kappa_p$
is given by \eqref{28-kappa}.  This is obvious when $r=0$. 
Next we suppose that $r=1$ and
$\kappa_p=0$. Then since $(k,Q)$ satisfies weak \textsf{LSC}, so
there is a solution to the congruence $a_1x_1^2\equiv k \tmod{p}$. Thus one
must have $(\frac{ka_1}{p})=1$, which again gives
an impossible  solution of \eqref{ganga}.
Finally, if $r=2$ and $\kappa_p=0$, then
Lemma~\ref{mr} implies that $M_2(p)>0$, 
which is again contrary to our assumption.
We may therefore assume that in any solution to the congruence 
$Q(\x)\equiv k \tmod{p}$, one has $p \mid \h(A_1x_1,\ldots,A_nx_n)$,
and there are at most two indices $i,j$ such that 
$p\nmid A_iA_j$ and $p\mid \h(x_i,x_j)$.

On returning to the diagonalisation process modulo $p$, we
deduce that there are linear forms $L_1,\ldots,L_r \in \Z[\X]$,
such that  $L_1(\x) \equiv \cdots \equiv L_r(\x)\equiv 0 \tmod p$ 
when $Q(\x) \equiv k \tmod p$.  Let
$$
\mathsf{\Lambda} = \{ \x \in \Z^n: p\mid L_1(\x),\ldots, p\mid L_r(\x)\},
$$ 
and recall that $0\leq r \leq 2$. Then we see that $\mathsf{\Lambda}$ forms a
lattice of determinant  $d(\mathsf{\Lambda})\in \{1,p,p^2\}$. In particular, there is a
basis $\mathbf{x}_1, \ldots, \mathbf{x}_n$ of $\mathsf{\Lambda}$ such that the matrix
$\ma{T}$ formed with $\mathbf{x}_1, \ldots, \mathbf{x}_n$ as column vectors has
$|\det \ma{T}| \in \{1,p,p^2\}$. We define 
$\widetilde{Q}(\y)=Q(\ma{T}\x)$, and note that by construction all coefficients of $\widetilde{Q}$
are divisible by
$p$. Thus we may write $\widetilde{Q}=pQ'$, where $Q'\in
\Z[X_1,\ldots,X_n]$ is a positive definite
quadratic form. Since $p\mid k$, we may also write $k=pk'$ for some
positive integer $k'$, and can proceed  to consider
the solubility of the equation $Q'(\x)=k'$ in integers.
Clearly $\mcal{S}(k;Q)=\emptyset$ if and only if $\mcal{S}(k';Q')=\emptyset$.
Furthermore it is trivial
to see that 
$$
|\Delta_{Q'}|=p^{-n}|\Delta_{\widetilde{Q}}|=
p^{-n}(\det \ma{T})^2 |\D| 
\leq |\D|,
$$
and $k/|\D|=p^{1-n} (\det \ma{T})^2 k'/|\Delta_{Q'}|\leq
k'/|\Delta_{Q'}|$,
since $n \geq 5$.
It remains to observe that $(k',Q')$ satisfies weak \textsf{LSC}, since $(k,Q)$ does.  
Indeed, in view of the fact that $\x \in \mathsf{\Lambda}$ whenever
$Q(\x) \equiv k \tmod p$, we see that 
$k'$ is represented by $Q'$ modulo $p^t$ for every $t \in \N$.  
%T%But this follows immediately from the observation that 
%T%$\x \in \mathsf{\Lambda}$ whenever $Q(\x) \equiv k \tmod p$, whence
%T%the transformation $\ma{T}$ was forced modulo $p^t$.

We now consider the problem of $2$-adic solubility. We shall show that
either $\sigma_2(k,Q)\gg 1$, or else we can replace the pair $(k,Q)$ by $(k',Q')$,
with $k'<k$, in such a way that $(k',Q')$ satisfies both weak \textsf{LSC}
and the conditions (ii) and (iii) of the lemma.
Quadratic forms are generally no longer 
diagonalizable over $\Z_2$.  However they are
``almost-diagonalizable'' (see Cassels \cite[Chapter 8, Lemma 4.1]{C}, for example),
in the sense that they can be written as a sum of forms of the type $aX^2$, $bXY$ and
$c(X^2+XY+Y^2)$, where $b$ and $c$ are even.
We may therefore proceed under the assumption that
\begin{equation}\lab{8-2adic}
\begin{split}
Q(\X) = & \sum_{i=1}^m a_i X_i^2 +
  \sum_{i=1}^{n_1}b_i X_{m+2i-1} X_{m+2i} \\
 &\quad +  \sum_{i=1}^{n_2} c_i\big(X_{m+2n_1+2i-1}^2 + X_{m+2n_1+2i-1} X_{m+2n_1+2i}
  + X_{m+2n_1+2i}^2\big),
\end{split}
\end{equation}
where $m+2n_1+2n_2=n$, and all $b_i$ and $c_i$ are even. We now
distinguish several cases. 

Suppose first that no coefficient $a_i$, $b_i$ or $c_i$ is divisible by
$8$.  We split our argument according to whether or not $2^7$ divides $k$.
Let us deal with the case $2^7\nmid k$ first.  By weak \textsf{LSC}
there must be a solution of $Q(\x) \equiv k \tmod{2^{13}}$ in which some $x_j$ is not
divisible by $2^4$.  Thus $2^7 \nmid \nabla Q(\x)$, by the
hypothesis that none of $a_i,b_i,c_i$ are divisible by $8$. Lifting this solution (see
\cite[Lemma 2.3]{dav32}, for example), we conclude that $\sigma_2(k,Q) \gg
1$.  Suppose next that $2^7\mid k$. Now clearly there is a solution
of $Q(\mathbf{x}) \equiv 0 \equiv k \tmod{2^7}$ where some $x_j$ is odd.
Indeed, since $n \ge 5$, the equation $Q=0$ is non-trivially soluble
in $\mathbb{Z}_2$, and thus also has a primitive solution. As above, we
conclude that $2^4 \nmid \nabla Q(\mathbf{x})$, and so deduce that
$\sigma_2(k,Q) \gg 1$ by a lifting argument.

We must now deal with the case in which at least one coefficient $a_i$, $b_i$ or $c_i$ is divisible
by $8$. There are two subcases to consider, the first being that
$Q(\x) \equiv k \tmod{2^7}$ has a solution
$\x$, such that there is an index $j \in \{1,\ldots,n\}$ for which both $2 \nmid x_j$
and $8 \nmid d_j$, where $d_j$ is the coefficient out of $a_i$, $b_i$,
$c_i$ that corresponds to $x_j$.  But here a lifting argument again
leads to the conclusion that $\sigma_2(k,Q) \gg 1$.
Finally, we consider the possibility that every solution of $Q(\x)
\equiv k \tmod{2^7}$ has $2\mid x_j$ or $8 \mid d_j$,  in the above notation, for each
$1\leq j\leq n$. In particular $k$ must be divisible by $4$.
Set $L_j(X_1, \ldots, X_n)$ equal to $0$ or $X_j$,
according to whether or not $8\mid d_j$, respectively.
Then the congruence $Q(\x) \equiv k \tmod{2^7}$ forces $2\mid L_j(\x)$ 
for $1\leq j \leq n$. Moreover, by our hypothesis at least one of
the $L_j$ must be identically zero. Hence the lattice
$$
\mathsf{\Gamma} = \{\x \in \mathbb{Z}^n: 2\mid L_1(\mathbf{x}),
\ldots, 2\mid L_n(\x)\},
$$
has determinant $d(\mathsf{\Gamma}) \le 2^{n-1}$. 
The pair $(k,Q)$ can now be reduced to $(k',Q')$ where
$k'=k/4$ and $Q'(\X)=4^{-1}Q(\ma{T}\X)$, with $|\det \ma{T}| \le 2^{n-1}$.
It follows that $|\Delta_{Q'}| = 4^{-n} |\det \ma{T}|^2 |\Delta_{Q}|
\le 4^{-1} |\D|$, whence
\[
  \frac{k'}{|\Delta_{Q'}|} \ge \frac{k/4}{|\D|/4} = \frac{k}{|\D|}.
\]
Furthermore,  $(k',Q')$ clearly satisfies weak \textsf{LSC}.  This
completes the proof of Lemma \ref{reduction}.
\end{proof}

Let $Q\in \Z[X_1,\ldots, X_n]$ be a positive definite quadratic form
and let $c_n>0$ be a fixed parameter that depends only upon $n$.
We proceed to introduce a quantity
$\k_n^\dagger(Q;c_n)$, that is defined to be the maximal $k\in\N$ such
that $\mcal{S}(k;Q)= \emptyset$, $\sigma_2(k,Q)\geq c_n$, and $(k,Q)$
satisfies strong $\textsf{LSC}$ modulo every odd prime $p$.  It is
not hard to see that the statement of Theorem \ref{w} remains valid
when $\k_n^*(Q)$ is replaced by $\k_n^{\dagger}(Q;c_n)$, with the
implied constant now depending additionally on the choice of $c_n$.
We claim that it will suffice to establish the upper bound in Theorem
\ref{iterate}, with $\k_n(Q)$ replaced by $\k_n^\dagger(Q;c_n)$, for a
suitable absolute constant $c_n>0$. 
To see this, let $k\in\K_n(Q)$. 
Then Lemma~\ref{reduction} implies that there is a pair $(k', Q')$
satisfying condition (i), such that $k'$ and $\Delta_{Q'}$ satisfy the
inequalities in part (iii) of the lemma, and the equation $Q'=k'$ has no solution
in integers. 
Suppose for the moment that we have shown $\k_n^\dagger(R;c_n)$ to be
$O(|\Delta_R|^{\varphi(n)})$, for $\varphi(n)\geq 1$ and
arbitrary positive definite quadratic forms $R \in
\Z[X_1,\ldots,X_n]$. Then we deduce that
$$
k\leq \frac{k'|\D|}{|\Delta_{Q'}|} \ll
|\Delta_{Q'}|^{\varphi(n)-1}|\D| \leq |\D|^{\varphi(n)}.
$$
This therefore establishes the claim.

Let us henceforth suppose that we have a positive integer $k$ such
that $k=\k_n^\dagger(Q;c_n)$, where $c_n>0$ is the absolute constant that
emerges from the application of Lemma \ref{reduction} above.
Our aim is to provide an upper bound for the size of $k$ in terms of the discriminant of $Q$.
It will suffice to replace $Q$ by any quadratic form $Q'$ that is
equivalent to it.  We may therefore proceed under the assumption that
$Q$ is Minkowski reduced.  Thus, if 
$$
Q(\X)=\sum_{1\leq i,j \leq n}q_{ij}X_iX_j,
$$
for appropriate integers $q_{ij}=q_{ji}$, then we may assume without loss of
generality that 
$$
\min_{\x\in\Z^n\setminus\{\ma{0}\}}Q(\x)=q_{11}\leq q_{22}\leq \cdots
\le q_{nn} \ll \Mq,
$$
with $q_{11}q_{22}\cdots q_{nn}\ll |\D|$.  Let us write 
$\min(Q)$ for the minimum non-zero value $Q$.  Then these inequalities
imply that
\beq\lab{min-max}
\min(Q)^{n-1}\Mq\ll |\D|.
\eeq
We now deduce from the
statement of Theorem~\ref{w},
with $\k_n^*(Q)$ replaced by $\k_n^{\dagger}(Q;c_n)$, that
\beq\lab{28-present}
k \ll_{\ve}
|\D|^{2/(n-3)}\Mq^{2n/(n-3)+\ve}.  
\eeq
This estimate is clearly at its sharpest when $\Mq$
is small compared to $|\D|$. It turns out that Watson's
approach produces a bound for $k$ that is best when $\Mq$ is large. 
Our plan is to combine the two bounds, in order to get an
overall improvement.
At this point it is convenient
to introduce a parameter $\al \in \R$, chosen so that 
\beq\lab{parameter}
A_1 |\D|^{\al} \leq \min(Q) \leq A_2 |\D|^{\al},
\eeq
for appropriate constants $A_2\geq A_1 \geq 0$ that depend only on $n$.
It is clear that this is always possible, and that $\al$ may be
taken to lie in the interval $[0,1/n]$.  With this convention it then
follows from \eqref{min-max} that $\Mq \ll |\D|^{1-\al(n-1)}$,
whence \eqref{28-present} yields
\beq\lab{f1}
k \ll_{\ve} |\D|^{2(n+1-\al n(n-1))/(n-3)+\ve},
\eeq
for any $\ve>0$.

We now turn to the bound for $k$ that emerges through an application
of Watson's method.  A cursory analysis of \cite[\S 9]{W} reveals the existence of positive
integers $a_1,\ldots,a_5$ such that $a_1=\min(Q)$ and 
\beq\lab{ho}
k\ll_{\ve} |\D|^{\max\{1,9/n\}}+a_1a_2a_3a_4a_5\big(|\D|^{\ve} + a_1^{-3/(n-3)}|\D|^{1/(n-3)}\big).
\eeq
Moreover it follows from \cite[Lemma 9.3]{W} that the integers
$a_1,\ldots,a_5$ satisfy the inequality
$$
a_2a_3a_4a_5 \ll (a_1|\D|)^{4/(n-4)}.
$$
On combining this with \eqref{parameter} and \eqref{ho}, we are
therefore led to the conclusion that
\begin{align*}
k
&\ll_{\ve} |\D|^{\max\{1,9/n\}}+
a_1^{n/(n-4)} |\D|^{4/(n-4)}\big(|\D|^\ve + 
a_1^{-3/(n-3)}|\D|^{1/(n-3)}\big)\\
&\ll |\D|^{\max\{1,9/n\}}+
|\D|^{(4+\al n)/(n-4)+(1-3\al)/(n-3)},
\end{align*}
provided that $\ve>0$ is chosen to be sufficiently small.
Let us assume that $n \leq 9$, since otherwise Watson's bound is
already best possible.  Then we deduce that
\beq\lab{f2}
k \ll |\D|^{9/n}+ |\D|^{(4+\al n)/(n-4)+(1-3\al)/(n-3)}.
\eeq

Define
$$
\al_0:=\frac{2n^2-11n+8}{2n^3-9n^2+2n+12}.
$$
It is not hard to check that \eqref{f1} is the sharper of the two
estimates for $\al\geq \al_0$, but that \eqref{f2} takes over for
$\al<\al_0$.  In this way we therefore deduce that 
\begin{align*}
k
&\ll_{\ve} 
|\D|^{9/n}+ |\D|^{(4+\al_0 n)/(n-4)+(1-3\al_0)/(n-3)+\ve}\\
&\ll_{\ve} |\D|^{9/n}+
|\D|^{\phi(n)+\ve},
\end{align*}
where $\phi(n)$ is given by \eqref{phi}.  This completes the
deduction of Theorem \ref{iterate}, since $\phi(n)\geq 9/n$ for
$5\leq n\leq 9$.  

\section{Activation of the circle method}\lab{act}

The purpose of this section is to establish Proposition
\ref{main'}. 
During the course of this we shall occasionally arrive at estimates
involving arbitrary parameters $M,N$.  These will typically be
non-negative or positive, but will always be assumed to take integer
values. 
We proceed to review the technical apparatus behind Heath-Brown's
version of the circle method \cite{HB'}.  Recall the
weight function $w_0:\R \rightarrow \R_{\geq 0}$, 
as given by \eqref{w0}, and set
$$
c_0:=\int_{-\infty}^\infty w_0(x) \d x.
$$
Let $\omega(x):=4c_0^{-1}w_0(4x-3)$, and define the function $h:
(0,\infty)\times \R \rightarrow \R$ by
$$
h(x,y):=\sum_{j=1}^\infty \frac{1}{xj}\Big(
\omega(xj)-\omega(|y|/xj)\Big).
$$
It is shown in \cite[\S 3]{HB'} that $h(x,y)$ is infinitely
differentiable for $(x,y) \in (0,\infty)\times \R$, and that $h(x,y)$
is non-zero only for $x \leq \max\{1,2|y|\}$.  
Let  $w_Q: \R^n \rightarrow \R_{\geq 0}$ be given by \eqref{wQ},
where $Q\in \Z[X_1,\ldots,X_n]$ 
is a non-singular quadratic form, as above.
The kernel of our work is \cite[Theorems 1 and  2]{HB'}.   For any $q \in \N$, and any $\c \in \Z^n$, we define the sum
\beq\lab{Sq} 
S_q(\c) := \sum_{\colt{a=1}{\hcf(a,q)=1}}^q
\sum_{\ma{b}\mod{q}}  e_q \big(a(Q(\ma{b})-k)+\ma{b}.\c\big), 
\eeq 
and the integral 
\beq\lab{Iq} 
I_q(\c):=
\int_{\R^n}w\Big(\frac{\x}{B}\Big)h\Big(\frac{q}{B},\frac{Q(\x)-k}{B^2}\Big)e_q(-\c.\x)\d\x. 
\eeq 
Then there exists a positive constant $c_B$, satisfying
$$
c_B=1+O_N(B^{-N}) 
$$
for any integer $N>0$, such that
\beq\lab{asym1} 
N_{w_Q}(Q-k;B) = c_B B^{-2}\sum_{\c\in \Z^n} \sum_{q
=1}^\infty q^{-n}S_q(\c)I_q(\c).  
\eeq

Our proof of Proposition \ref{main'} now has two major components: the estimation of the exponential sum
\eqref{Sq} and that of the integral \eqref{Iq}.  We shall treat these
separately, in \S \ref{sommes} and \S \ref{integrales}, respectively.
Finally, we shall deduce the statement of Proposition~\ref{main'} in
\S \ref{eichel}.

\subsection{Estimating $S_q(\c)$}\lab{sommes}

In this section we investigate 
the exponential sums $S_q(\mathbf{c})$, as given by
\eqref{Sq}.  We begin by recording the following basic
multiplicativity property \cite[Lemma 23]{HB'}.

\begin{lemma}\lab{mult}
If $\hcf(u,v)=1$ then  
$$ 
S_{uv}(\c)=S_{u}(\bar{v}\c) S_{v}(\bar{u}\c),  
$$
where $v \bar{v} \equiv {1} \tmod{u}$ and $u \bar{u} \equiv {1}
\tmod{v}.
$
\end{lemma}

The primary goal of this section is to obtain good upper bounds for $S_q(\c)$, in which
the implied constant is independent of the coefficients of $Q$. 
The following simple estimate is valid for
any choice of $q \in \N$.

\begin{lemma}
\label{expsumm}
We have 
$$
S_q(\mathbf{c}) \ll_{\ve}  q^{n/2+1+\ve}\hcf(q^n,\D)^{1/2}.
$$
\end{lemma}

\begin{proof}
We draw on the work of the second author \cite{dietmann}.  Let us write
$q=2^eq'$ and $\D=2^d \D'$, where $2\nmid q'\D'$ and $d,e$ are
non-negative integers. 
Then it follows from Lemma \ref{mult} that 
$$
S_q(\mathbf{c})=S_{2^e}(\bar{q'}\c)S_{q'}(\bar{2^e}\c)=S_1S_2,
$$
say.  Now \cite[Eq. (14)]{dietmann} immediately yields
$$
S_2\ll_{\ve} {q'}^{(n+1)/2+\ve}\hcf({q'}^n,\D')^{1/2}\hcf(q',k)^{1/2},
$$ 
whereas on combining \cite[Corollary 1]{dietmann} with the trivial estimate $|S_1| \leq 2^{en+e}$,
we see that
\begin{equation}\lab{31-S1}
S_1 \ll\min \big\{(2^e)^{n+1}, 2^{d/2}(2^e)^{n/2+1}\big\}
= (2^e)^{n/2+1}\hcf(2^{en}, 2^{d})^{1/2}.
\end{equation}
We therefore conclude that
$$
S_q(\c)\ll_{\ve} q^{n/2+1+\ve} \hcf({q'}^n,\D')^{1/2}\hcf(2^{en}, 2^{d})^{1/2}=q^{n/2+1+\ve}\hcf(q^n,\D)^{1/2},
$$
which thereby completes the proof of the lemma.
\end{proof}

We shall be able to achieve sharper bounds for $S_q(\mathbf{c})$
when $q$ is square-free.  Define the quadratic form 
$$
Q^{-1}(\x)=\x^T \mathbf{A}^{-1} \x, 
$$ 
with coefficients in $\Q$.  When $p$ is a prime such that $p \nmid
2\D$ we may think of $Q^{-1}(\x)$ as being defined modulo $p$. 
We now consider the sum $S_p(\c)$ for any odd prime $p$. 
By mimicking the argument
of \cite[Lemma 26]{HB'}, we  establish the following result.

\begin{lemma}\lab{s-p}
Let $p$ be an odd prime.
Then we have
$$
S_p(\c)\ll p^{(n+1)/2}\hcf(p^n,\D)^{1/2}\hcf(p,k, \del_n\D)^{1/2},
$$
where
\beq\lab{eta}
\del_n=\left\{\begin{array}{ll}
0, & \mbox{if $n$ even},\\
1, & \mbox{if $n$ odd}.
\end{array}
\right.
\eeq
\end{lemma}

\begin{proof}
Since $p$ is an odd prime there exists an integer valued matrix
$\mathbf{U}$ such that  $p \nmid \det \mathbf{U}$ and $\mathbf{U}^T\mathbf{A}\mathbf{U}$ is diagonal modulo $p.$  
In our estimation of $S_p(\c)$ it therefore suffices to assume that 
$Q(\X)=A_1X_1^2+\cdots+A_nX_n^2$, for integers $A_1,\ldots, A_n$ such that
$A_1\cdots A_n\equiv \D \tmod{p}$.
Suppose first that $p \nmid \D$.  In this setting it is not hard to see that
\begin{align*}
S_p(\c)
&=\sum_{a=1}^{p-1} e_p(-ak)\prod_{i=1}^n \sum_{b =1}^p
e_p( aA_ib^2+bc_i).
\end{align*}
One easily completes the treatment of the case $p\nmid \D$ by
recycling the arguments involving Gauss sums from the proof of Lemma
\ref{mr}, together with the well-known bounds
$$
|K_n(a,b;p)| \leq  \left\{
\begin{array}{ll}
2p^{1/2}\hcf(a,b,p)^{1/2}, & \mbox{if $n$ is even,}\\
p^{1/2}, & \mbox{if $n$ is odd,}\\
\end{array}
\right.
$$
Here, $K_n(a,b;p)$ is the Kloosterman sum for $n$
even, and the Sali\'e sum for $n$ odd.

Next we suppose that $p \mid \D$.  
On assuming that $Q$ has rank $\nu<n$ modulo $p$, the same sort of
argument leads to the conclusion that
$$
S_p(\c) \leq  p^{n-\nu}\Big| \sum_{a=1}^{p-1}
e_p(-ak) \prod_{i \in I} \sum_{b=1}^p e_p(aA_i b^2+ bc_i)\Big|,
$$
for some subset $I \subset [1,n]$ of cardinality $\nu$, such that $p
\nmid A_i$ for $i \in I$.  But then it is easy to deduce that
$$
S_p(\c) \ll p^{n-\nu}p^{(\nu+1)/2} \hcf(p,k)^{1/2}
\ll p^{(n+1)/2} \hcf(p^n,\D)^{1/2}\hcf(p,k)^{1/2}.
$$
Here, we have used the fact that
$p^{n-\nu}\mid \D$, since $F$ has rank $\nu$ modulo $p$, whence $n-\nu
\leq \min \{n,\nu_p(\D)\}$.  
This completes the proof of Lemma \ref{s-p}.
\end{proof}

We may now combine Lemma \ref{mult} and Lemma \ref{s-p} to provide an
estimate for $S_q(\mathbf{c})$ in the case that $q$ is square-free.

\begin{lemma}\label{expsumm-sfree}
Let $q \in \N$ be square-free. Then we have 
$$
S_q(\mathbf{c}) \ll_\ve
 q^{(n+1)/2+\ve} \hcf(q^n,\D)^{1/2}\hcf(q,k,\del_n \D)^{1/2},
$$
where $\del_n$ is given by \eqref{eta}.
\end{lemma}

\begin{proof}
Since $q$ is square-free we may write $q=2^e\prod_{j=1}^r p_j$, with
$p_1,\ldots,p_r$ distinct odd primes and $e \in \{0,1\}$.  Then it
follows from Lemma \ref{mult}, together with the trivial bound
$|S_2(\c)|\leq 2^{n+1}$, that
$$
|S_q(\ma{c})| \leq 2^{n+1} \prod_{j=1}^r |S_{p_j}(\bar{q_j}\ma{c})|.
$$
Here, $q_j=q/p_j$ for $1\leq j \leq r$, and $\bar{q_j}$ is defined by 
$q_j \bar{q_j} \equiv {1} \tmod{p_j}.$  
But then Lemma \ref{s-p} implies that 
$$
S_{p_j}(\bar{q_j}\c) \ll  
p_j^{(n+1)/2} \hcf(p_j^n,\D)^{1/2}\hcf(p_j,k, \del_n \D)^{1/2},
$$
for each $1 \leq j \leq r$.  On combining these two inequalities we
easily deduce the statement of Lemma \ref{expsumm-sfree}.
\end{proof}

We are now ready to investigate the average order of the sum $S_q(\c)$
for $q \leq X$.  To begin with we note that an application of Lemma
\ref{expsumm} immediately yields
\beq\lab{hoo}
\sum_{q \leq X}|S_q(\c)| \ll_{\ve} |\D|^{1/2} X^{n/2+2+\ve}.
\eeq
In fact we can do rather better than this in most circumstances.
Write $q=uv$ for coprime $u$ and $v$, such that
$u$ is square-free and $v$ is square-full.  
Then we may combine Lemmas \ref{mult}, \ref{expsumm} and
\ref{expsumm-sfree} to deduce that 
\begin{align*}
S_q(\ma{c}) &\ll_{\ve}  v^{n/2+1+\ve} |S_{u}(\bar{v}\c)| (v^n, \D)^{1/2}\\
&\ll_{\ve} u^{(n+1)/2}v^{n/2+1+\ve}
\hcf(u^nv^n,\D)^{1/2}\hcf(u,k,\del_n \D)^{1/2}\\
&\ll_{\ve} |\D|^{1/2}q^{(n+1)/2+\ve}v^{1/2} \hcf(u,k,\del_n \D)^{1/2},
\end{align*}
where $\del_n$ is given by \eqref{eta}.
Now for any non-zero integer $a$, and any $N \geq 1$, it is easy to see
that 
$$
\sum_{n\leq N}(n,a)\leq  \sum_{n \leq N}\sum_{e
\mid (n,a)}e  =\sum_{e \mid a}e\sum_{n'\leq N/e}1 \leq N d(a), 
$$
where $d(a)$ denotes the usual divisor function.
Let $C=\hcf(k,\del_n\D)$.  Then $C$ is a non-zero integer, unless $k=0$ and
$n$ is even.  Assuming this not to be the case, we employ the trivial
estimate $d(a)=O_\ve(|a|^\ve)$ in order to deduce that
\begin{align*}
\sum_{q \leq X}|S_q(\ma{c})| &\ll_{\ve}
|\D|^{1/2} X^{(n+1)/2+\ve}\sum_{v \leq X}  v^{1/2} \sum_{u \leq X/v} \hcf(u, C)\\
&\ll_{\ve}
|\D|^{1/2} C^\ve X^{(n+3)/2+\ve}\sum_{v \leq X}v^{-1/2}.
\end{align*}
On noting that there are only $O(V^{1/2})$ square-full values of $v
\leq V$, we have therefore established that
$$
\sum_{q \leq X}|S_q(\c)| \ll_{\ve} 
|\D|^{1/2+\ve}(1+k)^\ve X^{(n+3)/2+\ve},
$$
unless $k=0$ and $n$ is even. We may now combine this with
\eqref{hoo} in order to deduce the following result.

\begin{lemma}\lab{exp-av-q}
For any $X \geq 1$ we have
$$
\sum_{q \leq X}|S_q(\c)| \ll_{\ve} 
|\D|^{1/2+\ve}(1+k)^\ve X^{(n+3+\gamma_n)/2+\ve},
$$
where $\gamma_n$ is given by \eqref{even-odd}.
\end{lemma}

We end this section by considering the average order of the sum 
$S_q(\c)$ in the special case $\c=\ma{0}$.
But in this setting Lemma \ref{exp-av-q} clearly yields
$$
\sum_{q \leq X} q^{-n}S_q(\ma{0}) = \sum_{q=1}^\infty
q^{-n}S_q(\ma{0}) +O_{\ve}\big(|\D|^{1/2+\ve}(1+k)^\ve X^{(3+\gamma_n-n)/2+\ve}\big),
$$
which implies that the infinite sum in this formula is absolutely convergent.  
Lemma \ref{mult} implies that the function $q^{-n}S_q(\ma{0})$ is
multiplicative.  Thus the usual analysis of the singular series yields
\beq\lab{yemuna}
\sum_{q=1}^\infty
q^{-n}S_q(\ma{0})=\prod_p \sum_{t=0}^\infty p^{-nt}S_{p^t}(\ma{0})\\
=\prod_p \sigma_p,
\eeq
where $\sigma_p$ is given by \eqref{sig-p}.  On recalling the
definition \eqref{singseries} of $\ss(k,Q)$, we  have therefore shown that
\beq\lab{singseries1}
\sum_{q \leq X} q^{-n}S_q(\ma{0}) = \ss(k,Q) +O_{\ve}\big(|\D|^{1/2+\ve}(1+k)^\ve X^{(3+\gamma_n-n)/2+\ve}\big),
\eeq
where $\gamma_n$ is given by \eqref{even-odd}.

\subsection{Estimating $I_q(\c)$}\lab{integrales}

The goal of this section is to provide good upper bounds for the
integral $I_q(\c)$, for given $q \in \N$ and $\c \in \Z^n$. Recall the orthogonal matrix 
$\mathbf{R} \in \mathrm{SO}_n(\R)$ that was fixed at the outset.  Thus
there exist non-zero $\lambda_1,\ldots,\lambda_n \in \R$ 
such that  \eqref{choice-R} holds, with $\lambda_1\cdots \lambda_n=\D$.  It
will be convenient to introduce the matrix
$$
\mathbf{D}:=\mathrm{Diag}(\lambda_1^{-1/2},\ldots,\lambda_n^{-1/2}).
$$
Recall now the definitions \eqref{w1}, \eqref{wQ} and \eqref{tW} of the weight
functions $w_1, w_Q$ and $\tilde{w}$, respectively.  Finally, recall the
definition  \eqref{Iq} of the integral $I_q(\c)$, and also that of the
polynomial $P_Q$, given by \eqref{P}.  Then it follows from a simple change
of variables that
\begin{align*}
I_q(\c)
& = B^n\int_{\R^n}w_Q(\x)h\big(B^{-1}q,B^{-2}(Q(B\x)-k)\big)e_q(-B\c.\x)\d\x\\
& = B^n\int_{\R^n}\tilde{w}(\mathbf{R}^T\x)h(B^{-1}q,P_Q(\x))e_q(-B\c.\x)\d\x.
\end{align*}
But then the change of variables $\u=\mathbf{R}^T\x$ easily yields
\begin{align*}
I_q(\c)
& = B^n \int_{\R^n}
\tilde{w}(\u)h(B^{-1}q,P_Q(\mathbf{R}\u))e_q(-B(\mathbf{R}^T \c).\u)\d\u\\
& = \frac{B^n}{|\D|^{1/2}} \int_{\R^n}
w_1(\u) h(B^{-1}q,P_{Q_{\mathrm{sgn}}}(\u))e_q(-B\ma{v}.\u)\d\u,
\end{align*}
where $\ma{v}=\mathbf{D}\mathbf{R}^T \c$
and $Q_{\mathrm{sgn}}$ is given by \eqref{R-sig}.   In particular we have $\|P_{Q_{\mathrm{sgn}}}\|\leq 1$, so that 
$|P_{Q_{\mathrm{sgn}}}(\u)|\ll 1$ for any $\u \in \supp(w_1)$.
Once taken together with the properties of the function $h$ mentioned
in \S \ref{2}, we deduce that $I_q(\c)=0$ unless $q\ll B$.
Moreover, in view of the fact that $\mathbf{R}$ is orthogonal, it easily follows from \eqref{eigen} that
\beq\lab{jonny}
|\ma{v}|\gg |\c|/\Mq^{1/2},
\eeq
where $\Mq$ denotes the height of $Q$, as usual.

Now let $w: \R^n \rightarrow \R_{\geq 0}$ be any infinitely
differentiable function, with compact support, and let $G\in \R[\x]$
be any quadratic form.  Much of this section will be based on an
analysis of the integral
\beq\lab{I*}
I_r^*(\ma{v};G,w):= \int_{\R^n}w(\x)h(r,P_G(\x))e_r(-\ma{v}.\x)\d\x,
\eeq
for any $r \in \R$ and $\ma{v}\in \R^n$.  Here, $P_G$ is given by
\eqref{P} as usual, and so takes the form $P_G(\x)=G(\x)$ or 
$P_G(\x)=G(\x)-1$, according to whether we are interested in the case
$k=0$ or $k>0$ in our proof of Proposition \ref{main'}.
The connection with $I_q(\c)$ is given by 
\beq\lab{I^*}
I_q(\c)= \frac{B^n}{|\D|^{1/2}} I_{r}^*(\ma{v};Q_{\mathrm{sgn}},w_1),
\eeq
with $r=B^{-1}q$ and $\ma{v}=\mathbf{DR}^T \c$. 
In order to establish upper bounds for $I_q(\c)$, it will
therefore suffice to do so for  $I_r^*(\ma{v};G,w)$, 
for any quadratic form $G\in \R[\x]$, and a rather general class
$\CC_1(S)$ of functions $w$, that depend only upon parameters from a
set $S$.

The class of functions $\CC_1(S)$ that we shall work with is very
similar to that employed by Heath-Brown \cite[\S\S 2,6]{HB'}.
We shall use weights $w: \R^n \rightarrow
\R_{\geq 0}$, which are infinitely differentiable and have compact
support, and which take non-negative real values.  Given such a weight 
$w$, we let
$n(w)=n$ and set $\rom{Rad}(w)$ to be the smallest $R$
such that $w$ is supported in the hypercube $[-R,R]^n$.  Moreover for
every integer $j \geq 0$ we let
$$
\kappa_j(w):=\max\Big\{ \Big|
\frac{\partial^{j_1+\cdots+j_n}w(\x)}{\partial^{j_1}x_1\cdots
\partial^{j_n}x_n}\Big|: ~\x \in \R^n, ~j_1+\cdots+j_n=j\Big\}.
$$
Let $S$ be any collection of parameters.  Then we shall define $\CC(S)$ to be
the set of infinitely differentiable functions $w: \R^n \rightarrow
\R_{\geq 0}$ of compact support, such that $n(w),\rom{Rad}(w),
\kappa_0(w), \kappa_1(w), \ldots$ are all bounded by corresponding
quantities $n(S),\rom{Rad}(S), \kappa_0(S), \kappa_1(S), \ldots$
depending only on parameters from the set $S$.  
This much is in complete accordance with Heath-Brown  \cite[\S 2]{HB'}.
Given a quadratic form $G \in \R[\X]$,  we now specify the set of
functions $\CC_1(S)\subset \CC(S)$, where we now assume that $S$
contains among its parameters the form $G$.  Given $w \in \CC(S)$, we shall say that $w \in \CC_1(S)$ if 
$$
\frac{\partial G}{\partial X_1} \gg_{S} 1
$$
on $\supp(w)$.  In particular, when $G=Q_{\mathrm{sgn}}$ is given by \eqref{R-sig}, it is not hard to
see that $w_1 \in \CC_1(Q_{\mathrm{sgn}},n)=\CC_1(n)$.

We are now ready to commence our study of the integral
$I_r^*(\ma{v};G,w)$, for an arbitrary quadratic form $G\in \R[\X]$ and
a general weight function $w\in \CC_1(S)$, where as above 
$S$ is assumed to contain $G$ among its parameters.
Phrasing things in this degree of generality  allows us to 
apply the work of Heath-Brown more or less
directly, since we have a direct correspondence between \eqref{I*} and the integral
$I_r^*(\ma{v})$ in \cite[\S 7]{HB'}, 
defined for any quadratic form $G$ and weight $w\in \CC_1(S)$.
There is a slight abuse of notation here, in that the polynomial $G$
appearing in Heath-Brown's definition of $I_r^*(\ma{v})$ corresponds
precisely to what we have called $P_G$, for a quadratic form $G$.
Moreover, it should be highlighted that whereas his work is phrased in terms of the more restrictive
class of weight functions $\CC_0(S)\subset \CC_1(S)$, defined 
at the start of \cite[\S 6]{HB'}, an inspection of the contents of
\cite[\S\S 7,8]{HB'} reveals that all of the estimates there
extend to $\CC_1(S)$.  

Bearing this in mind, our first task is to 
record a preliminary estimate for $I_q(\c)$.  The following will be used to show that large values of $\c$ 
make a negligible contribution in our analysis.

\begin{lemma}\lab{c>0:1}
Let $\ma{c} \in \Z^n$ with $\c \neq \ma{0}$.  Then 
for any $N\geq 0$ we have 
$$
I_q(\c)\ll_{N}  \frac{B^{n+1}}{q|\D|^{1/2}} \frac{\Mq^{N/2}}{|\c|^{N}}.
$$
\end{lemma}

\begin{proof}
Let $G\in \R[\X]$ be a quadratic form and let $w\in \CC_1(S)$.  Then on combining
\cite[Lemma 14]{HB'} and \cite[Lemma 18]{HB'}, we deduce that
$$
I_r^*(\ma{v};G,w) \ll_{N,S} r^{-1}|\ma{v}|^{-N},
$$
for any $N\geq 0$ and any $\ma{v} \in \R^n$ such that $\ma{v} \neq
\ma{0}$.   But then we may insert this into \eqref{I^*}, and combine
it with \eqref{jonny},  in order to complete the proof.
\end{proof}

We shall need a finer estimate for $I_q(\c)$ when $\c$ has small
modulus.  The following result also follows rather easily from
Heath-Brown's analysis.

\begin{lemma}\lab{c>0:2}
Let $\c \in \Z^n$ with $\c\neq \ma{0}$.
Then we have
$$
I_q(\c) \ll_{\ve} 
\frac{\Mq^{n/4-1/2+\ve}}{|\D|^{1/2}}
B^{n/2+1+\ve} q^{n/2-1}|\c|^{1-n/2+\ve}.
$$
\end{lemma}

\begin{proof}
For an arbitrary quadratic form $G\in \R[\X]$, the Hessian condition
of \cite[Lemma 21]{HB'} automatically holds for any $w\in
\CC_1(S)$.  Thus for any such $G$ and $w$, we may
combine \cite[Lemma 14]{HB'} with \cite[Lemma 22]{HB'} to deduce that
$$
I_r^*(\ma{v};G,w)  \ll_{\ve, S} (r^{-2}|\ma{v}|)^\ve r^{n/2-1} |\ma{v}|^{1-n/2},
$$
for any $\ma{v} \in \R^n$ such that $\ma{v} \neq
\ma{0}$.   On inserting this into \eqref{I^*}, and
combining it with \eqref{jonny},  we therefore complete the proof of
Lemma \ref{c>0:2}.
\end{proof}

We end this section by considering 
$I_q(\ma{0})=|\D|^{-1/2}B^nI_r^*(\ma{0};Q_{\mathrm{sgn}},w_1)$ in \eqref{I^*}, with
$r=B^{-1}q$.  
Let $G \in \R[\X]$ be an arbitrary quadratic form, and let $P_G \in
\R[\X]$ be the corresponding polynomial \eqref{P}.
Then it follows from \cite[Lemma~13]{HB'} that
$$
I_r^*(\ma{0};G,w)= \sigma_\infty(w;P_G) +O_{N,S}(r^N),
$$
for any $N>0$, where $\sigma_\infty(w;P_G)$ is given by \eqref{singint-1},
and we have assumed that $r \ll_S 1$.  On inserting this into
\eqref{I^*}, and recalling the identity \eqref{sig-inf},  we therefore obtain the following result.

\begin{lemma}\lab{c=0}
We have
$$
I_q(\ma{0})= \frac{\sigma_\infty B^n}{|\D|^{1/2}}  +O_{N}\Big(\frac{q^NB^{n-N}}{|\D|^{1/2}}\Big),
$$
for any $N>0$.
\end{lemma}

\subsection{Derivation of Proposition \ref{main'}} \label{eichel}

In this section we are going to complete the proof of Proposition \ref{main'}. 
Our starting point is \eqref{asym1}. Let $\ve>0$ and let  $P\geq 1$.  
Then, whether we are in the case $k=0$ or $k>0$, we always have 
$(1+k)^\ve\ll B^\ve$ in Lemma \ref{exp-av-q}. On combining this result with Lemma \ref{c>0:1},  and the fact that
$I_q(\c)=0$ unless $q\ll B$, we therefore see that the contribution to the right hand
side of \eqref{asym1} from  $|\c|> P$ is 
\begin{align*}
&\ll_{N} 
\frac{B^{n-1}}{|\D|^{1/2}} 
\sum_{|\c|> P}
\sum_{q\ll B} q^{-n-1}|S_q(\c)|  
  \frac{\Mq^{N/2}}{|\c|^{N}}
\ll_{\ve,N}
B^{n-1+\ve} \frac{\Mq^{N/2+\ve}}{P^{N-n}},
\end{align*}
for any $N> n$.  But this is clearly 
\beq\lab{fri=1}
\ll_{\ve, M}
B^{n-1+\ve} \frac{\Mq^{(M+n)/2+\ve}}{P^{M}},
\eeq
for any $M>0$.
Turning to the contribution from $1\leq |\c|\leq P$,
we combine Lemmas \ref{exp-av-q} and \ref{c>0:2} to deduce that 
\begin{align*}
\sum_{q=1}^\infty q^{-n}S_q(\c)I_q(\c)
&\ll
\max_{Y \ll B}\sum_{j\leq \log Y}
\sum_{2^{j-1}< q \leq 2^j}q^{-n}|S_q(\c)I_q(\c)| \\
&\ll_{\ve} 
\Mq^{n/4-1/2+\ve} B^{(n+3+\gamma_n)/2+\ve} |\c|^{1-n/2+\ve},
\end{align*}
where $\gamma_n$ is given by \eqref{even-odd}.
On summing over values of $\c$ such that $1\leq |\c|\leq P$,
we therefore deduce that
the contribution to the right hand side of \eqref{asym1} from
such $\c$ is 
$$
\ll_{\ve} \Mq^{n/4-1/2+\ve} B^{(n-1+\gamma_n)/2+\ve}P^{n/2+1+\ve}.
$$
Once combined with \eqref{fri=1}, we see that
the overall contribution from $\c\neq \ma{0}$ is 
$$
\ll_{\ve,M}
\Mq^{n/4-1/2+\ve}B^{(n-1+\gamma_n)/2+\ve}\Big(
\frac{\Mq^{M/2+n/4+1/2}B^{n/2}}{P^{M}} + 
P^{n/2+1+\ve}\Big),
$$
for any $M> 0$.  Taking $M=\lceil n/(2\ve)\rceil$ and 
$P= \Mq^{1/2} B^{\ve},$ 
we therefore see that there is a contribution of
\beq\lab{water}
\ll_{\ve}
\Mq^{n/2+\ve}B^{(n-1+\gamma_n)/2+\ve}
\eeq
to the right hand side of \eqref{asym1} from
those $\c\neq \ma{0}$.

It remains to handle the contribution from the case $\c=\ma{0}$.
Recall the definition \eqref{sig-inf} of $\sigma_\infty$, and
the inequality $\sigma_\infty\ll 1$ that it satisfies.
Then an application of Lemma \ref{c=0}, together with  \eqref{singseries1}, reveals that the contribution 
from $\c= \ma{0}$  is 
\begin{align*}
\frac{c_B}{B^2}\sum_{q \ll B} q^{-n}S_q(\ma{0})I_q(\ma{0})
=&
\frac{\sigma_\infty \ss(k,Q) B^{n-2}}{|\D|^{1/2}} +O_{\ve}(B^{(n-1+\gamma_n)/2+\ve})\\
&\quad + O_{N}\Big(\frac{B^{n-2-N}}{|\D|^{1/2}} \sum_{q \ll B} q^{N-n}|S_q(\ma{0})|\Big),
\end{align*}
for any $N>0$.    On selecting $N=(n-3-\gamma_n)/2$, it follows
from Lemma \ref{exp-av-q} that the error terms in this estimate are
bounded by $O_{\ve}(B^{(n-1+\gamma_n)/2+\ve}).$
We may now combine this with \eqref{water} in \eqref{asym1}, in order to complete the proof of 
Proposition \ref{main'}.

\section{The singular series}\label{ss}

In this section we establish Proposition \ref{ss'}.  
Let $Q \in \Z[X_1,\ldots,X_n]$ be a non-singular quadratic form
of discriminant $\D$, and let $k$ be a non-negative integer. As usual
we assume that $n\geq 5$ when $k=0$. Under suitable local solubility
assumptions,  our task is to establish a uniform lower bound for the singular
series $\ss(k,Q):=\prod_p \sigma_p$, where 
$\sigma_p=\sigma_p(k,Q)$ is given by \eqref{sig-p} and \eqref{N}.
We shall obtain a better lower bound by assuming stronger local
solubility conditions.
Note that when $n\geq 5$ it already follows that
$\sigma_p \ne 0$ under weak \textsf{LSC}.
Similarly, when $n=4$ and $k>0$, it follows from 
strong \textsf{LSC} that $\sigma_p \ne 0$.

Let us begin by handling the factors $\sigma_p$, for which $p
\nmid \D$.  In this setting we shall use the identity
$$
\sigma_p=\sum_{t=0}^\infty p^{-nt}S_{p^t}(\ma{0}),
$$
that follows from \eqref{yemuna}.
But then an application of Lemma \ref{expsumm} reveals that 
$$
\sigma_p = 1+O_{\ve}\Big(\sum_{t\geq 1}p^{t(1-n/2+\ve)} \Big)=1+O_{\ve}(p^{-3/2+\ve}),
$$
when $n \geq 5.$  When $n=4$ and $k>0$, an application of
Lemmas \ref{expsumm} and \ref{s-p} also yields
$\sigma_p = 1+O_{\ve}(p^{-3/2+\ve})$, provided  that $p\nmid 2k\D$.
Hence we have
\beq\lab{mendip1}
\prod_{p \nmid \D} \sigma_p \gg 1, \quad (n\geq 5),
\eeq
and 
\beq\lab{mendip2}
\prod_{p \nmid 2k\D} \sigma_p \gg 1, \quad (\mbox{$n=4$ and  $k>0$}).
\eeq

We now turn to the size of $\sigma_p$ for the remaining primes $p$.
For this we recall the definitions \eqref{N} and \eqref{N*} of
$N(p^t)$ and $N^*(p^t)$, respectively.
When $p$ is odd, we may diagonalise $Q$ modulo
$p^t$, without changing the values of $N(p^t)$ or $N^*(p^t)$.  In our analysis of
the quantities $N(p^t)$ or $N^*(p^t)$ for odd $p$, it therefore suffices to proceed 
under the assumption that 
\beq\lab{7-diag}
Q(\X)=A_1X_1^2+\cdots+A_nX_n^2,
\eeq
with $A_1\cdots A_n \equiv \D \tmod{p^t}.$  
Now it is easy to deduce from \eqref{sig-p} and \eqref{hensel} that in the case $n=4$
and $k>0$ we have
\beq\lab{new-house}
\sigma_p\geq p^{-3(1+2\tau_p)}N^*(p^{1+2\tau_p}).
\eeq
When $p=2$, therefore, we immediately obtain $\sigma_2\geq 2^{-9}$, 
since the pair $(k,Q)$ is assumed to satisfy strong \textsf{LSC}
modulo $2$.  Suppose next that $p>2$, with $p\mid k\D$.
We may therefore assume that $Q$ takes the shape \eqref{7-diag},
with $A_1A_2A_3A_4\equiv \D \tmod{p}$.
On combining  Lemma \ref{mr} with the fact that $(k,Q)$ satisfies strong
\textsf{LSC}, one may deduce that 
$N^*(p)=2p^3$ if $p\mid \hcf(k,A_i,A_j)$ for precisely two indices
$1\leq i<j\leq 4$, and $N^*(p)=p^3+O(p^2)$ otherwise.
We shall revisit this line of argument in greater detail when we deal
with the case $n\geq 5$ below.
On inserting this into \eqref{new-house} and combining it with
\eqref{mendip2} and our lower bound for $\sigma_2$, we therefore deduce that
there is an absolute constant $c>0$ such that
$$
\ss(k,Q)\gg \prod_{p\mid k\D}
\Big(1-\frac{c}{p}\Big)
\gg_\ve k^{-\ve}|\D|^{-\ve},
$$
for any $\ve>0$.
This is satisfactory for the statement of Proposition \ref{ss'} when 
$n=4$ and $k>0$.

We may assume henceforth that $n\geq 5$.  Suppose first that $p$ is
odd with $p\mid \D$. Then we may assume that $Q$ takes the shape
\eqref{7-diag}, with
$
A_1\cdots A_n \equiv \D \tmod{p^t}.
$  
We shall say that the pair $(k,Q)$ is ``$p$-reduced'' if any of the
following occur:
\begin{enumerate}
\item $p$ divides at most $n-3$ of the coefficients $A_1,\ldots,A_n$.
\item $p$ divides $k$ and precisely $n-2$ of the coefficients
$A_1,\ldots,A_n$, with the remaining two coefficients satisfying 
$(\frac{-A_iA_j}{p})=1$.
\item $p$ divides all but one of the coefficients
$A_1,\ldots,A_n$, with the remaining coefficient satisfying 
$(\frac{kA_i}{p})=1$.
\end{enumerate}
We shall be able to establish a good lower bound for $N(p^t)$ when 
the pair $(k,Q)$ is $p$-reduced.  When $(k,Q)$ is not
$p$-reduced, we will be able to make a certain change of variables that ultimately
leads us to estimate $N(p^t)$ for a $p$-reduced pair $(k',Q')$.
Observe that any pair $(k,Q)$ satisfying strong \textsf{LSC}, is
automatically $p$-reduced, since then $N^*(p)>0$.
Suppose for the moment that $(k,Q)$ satisfies weak \textsf{LSC}, but
is not $p$-reduced.  
We claim that there exists an integer
$\theta\geq 1$ and a $p$-reduced pair $(k',Q')$ such that
\beq\lab{31-pred}
N(p^t)\geq p^{\theta(n-2)} \#\{\x \mod{p^{t-\theta}}: Q'(\x)\equiv k' \mod{p^{t-\theta}}\},
\eeq
with
\beq\lab{31-theta}
\theta\leq \frac{\nu_p(\D)}{n-4}.
\eeq
To see this we note that after a possible relabelling of the indices
we have
$$
(A_1,\ldots,A_n)=(a_1,\ldots,a_r,p^{\be_1}b_1,\ldots,p^{\be_s}b_s),
$$
where $(r,s)=(0,n),(1,n-1)$ or $(2,n-2)$, and 
\beq\lab{1-aibi}
p \nmid a_ib_i, \quad  1\leq \be_1\leq \cdots \leq \be_s, \quad
\be_1+\cdots+\be_s=\nu_p(\D).
\eeq
If $r=0$ then we may write $k=pk'$ in $N(p^t)$, since $N(p)>0$ by assumption.  This gives
$$
N(p^t)
=p^{n}\#\{\x \mod{p^{t-1}}: ~Q_0(\x) \equiv k' \mod{p^{t-1}}\},
$$
where $Q_0(\ma{x})=p^{\beta_1-1}b_1X_1^2+\cdots+p^{\beta_{n}-1}b_{n}X_{n}^2$.
If $r=1$ and $(\frac{ka_1}{p})\neq 1$, then we may write $x_1=px_1'$
and $k=pk'$ in $N(p^t)$, again since $N(p)>0$.  This gives
\begin{align*}
N(p^t)
&=p^{n-1}\#\{\y \mod{p^{t-1}}: ~Q_1(\y) \equiv k' \mod{p^{t-1}}\},
\end{align*}
where $Q_1(\ma{Y})=p^{\beta_1-1}b_1Y_1^2+\cdots+
p^{\beta_{n-1}-1}b_{n-1}Y_{n-1}^2+pa_1{Y_n}^2$.
Alternatively, if $r=2$ and $(\frac{-a_1a_2}{p})= -1$, then we may write
$x_1=px_1', x_2=px_2'$
and $k=pk'$ in $N(p^t)$.  This gives
\begin{align*}
N(p^t)
&=p^{n-2}\#\{\z \mod{p^{t-1}}: ~Q_2(\z) \equiv k' \mod{p^{t-1}}\},
\end{align*}
where $Q_2(\ma{Z})=p^{\beta_1-1}b_1Z_1^2+\cdots+
p^{\beta_{n-2}-1}b_{n-2}Z_{n-2}^2+pa_1{Z_{n-1}}^2+pa_2{Z_{n}}^2$.
In particular it is clear that 
$$
\nu_p(\Delta_{Q_i})= \nu_p(\D)-n+2i, \quad
(\mbox{$i=0,1,2$}).
$$
Now either $(k',Q_i)$ is $p$-reduced, or we can repeat the
procedure. This procedure must eventually terminate, since at each step the
$p$-adic valuation of the forms determinant is reduced by $\geq
n-4$, and we can stop when the forms determinant has $p$-adic
valuation $\leq n-3$. More precisely, after the $j$th step one finds that
the $p$-adic valuation of the resulting form is $\leq
\nu_p(\D)-j(n-4)$, whence there at most 
$$
1+\Big[\frac{\nu_p(\D)-n+3}{n-4}\Big]
$$
steps in this
procedure. The bounds in \eqref{31-pred} and \eqref{31-theta} are now obvious.
We note that in the alternative case, where $(k,Q)$ satisfies strong \textsf{LSC}, 
we take $\theta=0$ and $(k',Q')=(k,Q)$ in \eqref{31-pred}.

Still under the assumption that $p>2$, with $n\geq 5$ and $p\mid \D$,
we proceed to derive a lower bound for 
$N(p^t)$ for any $t\in \N$, under the assumption that
the pair $(k,Q)$ is $p$-reduced.  Moreover we shall suppose that
$$
Q(\X)=a_1X_1^2+\cdots + a_rX_r^2+p^{\be_1}b_1X_{r+1}^2+\cdots+p^{\be_s}b_sX_n^2,
$$
with $r,s\in\N$ such that $r+s=n$, and \eqref{1-aibi} holding.
%T%include?
%% \begin{proof}
%% Let  $\x \tmod {p^{u-1}}$ be such that $Q(\x)-k=p^{u-1} a$, say, and write
%% $
%% \x'=\x+p^{u-1-\tau_p}\ma{y}
%% $
%% for some $\ma{y} \tmod {p}$.
%% Then in order to obtain a lower bound for $N'(p^{u})$, it suffices to
%% find values of $\ma{y}$ for which 
%% $$
%% 0 \equiv Q(\x')-k \equiv p^{u-1}a+ p^{u-1}\nabla Q(\x).\ma{y} \mod{ p^{u}}.
%% $$
%% Noting that $\nabla Q(\x).\ma{y}=2 (\ma{Ax}).\y$, it is not hard to
%% see that  there are at least $p^{n-1}$ possible such
%% $\ma{y} \tmod{p}$.  This therefore establishes
%% (\ref{hensel}). 
%% \end{proof}
It follows from \eqref{hensel} that
\beq\lab{1-food}
N(p^t)\geq p^{(n-1)(t-1)}N^*(p),
\eeq
for any $t > 1$.  Observe that
\begin{align*}
N^*(p)
&= p^{s}\#\{\z \mod{p}:
~a_1z_1^2+\cdots+a_rz_r^2\equiv k \mod{p}, ~p\nmid \z\}.
\end{align*}
Hence $N^*(p)=p^sM_r(p)$, in the notation of \eqref{mr-def}.
By combining this fact 
with Lemma \ref{mr} we deduce that
\beq\lab{7-lower}
N^*(p) \geq p^{n-1}+O(p^{n-2}).
\eeq
Let $t$ be large and suppose that $(k,Q)$ satisfies weak \textsf{LSC}.  
Then on combining \eqref{31-pred} and \eqref{1-food} with \eqref{7-lower}, we may conclude that
$$
N(p^t) \ge p^{t(n-1)-\theta}(1+O(p^{-1})),
$$
for some $\theta\geq 1$ such that \eqref{31-theta} holds.
Thus we deduce that 
$$
\sigma_p=\lim_{t\rightarrow \infty} p^{-t(n-1)}N(p^t)\geq p^{-\theta}(1+O(p^{-1}))\geq p^{-\nu_p(\D)/(n-4)}(1+O(p^{-1}))
$$
in this case.  Alternatively, if 
$(k,Q)$ satisfies strong \textsf{LSC}, so that one may take $\theta=0$
in \eqref{31-pred}, it follows from \eqref{1-food} and \eqref{7-lower} that 
$$
\sigma_p\geq 1-\frac{c_n}{p},
$$
for some constant $c_n>0$ depending at most upon $n$.

We now turn to the case $p=2$, which will be handled by an iterative
method similar to that used above. If $k$ is not divisible by $32$
then we claim that 
\beq\lab{N2}
N(2^t) \gg 2^{t(n-1)}.
\eeq
This clearly implies that $\sigma_2\gg 1$. To see the claim
we note that $\nabla Q(\mathbf{x})$ cannot be divisible by $16$ for
any $2$-adic solution of $Q(\mathbf{x})=k$, and such a $2$-adic
solution must exist by weak \textsf{LSC}. 
A lifting argument therefore establishes the claim.
We proceed under the assumption that $k$ is divisible by $32$. 
We may assume that $Q$ takes the shape \eqref{8-2adic} for certain
integers $a_i,b_i,c_i$ such that $b_i$ and $c_i$ are even.  Suppose
that there is at least one coefficient among the $b_i,c_i$ that is not
divisible by $8$, and write $F(X,Y)$ for the binary quadratic form that
corresponds to this coefficient.
Then it is not hard to see that there is a solution of 
$$
F(x,y) \equiv 0 \equiv k \mod {32},
$$ 
with $2\nmid x$.  Thus we may find a solution of $Q(\x)\equiv
k\tmod{32}$ in which $\nabla Q(\mathbf{x})$ is not divisible by
$8$. This solution can be lifted, which thereby shows that \eqref{N2} holds in this
case also. We may henceforth assume that every coefficient $b_i,c_i$
in \eqref{8-2adic} is divisible by $8$. In particular they may be divided by $4$
and still leave terms that are classically integral.

We now repeat the argument that we used to treat the case
$p> 2$, although extra care needs to be taken to obtain a result of the
same strength. If in the course of the reduction only one or two
variables are forced to be even in order that
$Q(\mathbf{x}) \equiv 0 \equiv k \tmod{2}$, then we may proceed
analogously to $p>2$, obtaining $(k',Q')$ in place of $(k,Q)$, such that
\eqref{31-pred} holds with $p=2$ and $\theta=1$.
We then enter the next iteration step, repeating our argument for
$p>2$, but now with $(k',Q')$.  In the alternative case, the
congruence $Q(\mathbf{x}) \equiv 0 \tmod{2}$
forces more than two variables to be even. But then our quadratic form
$Q$ splits off a diagonal form 
$$
a_1 X_1^2 + \cdots +a_m X_m^2,
$$ 
where $m \ge 3$ and $2\nmid a_1\cdots a_m$, and $x_1, \ldots, x_m$ are forced to
be even in any solution to the congruence $Q(\mathbf{x}) \equiv 0 \tmod{2}$.
Note that no non-diagonal terms appear here, by our
preparations above. Clearly we may assume that 
$m \le 4$, since for $m \ge 5$
all $k$ could be represented by $x_1, \ldots, x_m$ not all even. Thus
another lifting argument would yield \eqref{N2}.
We shall also need to note that if any
other coefficients associated to $X_{m+1}, \ldots, X_n$ were divisible
by $2$, but not by $4$, then we could solve the congruence
$$
Q(\x) \equiv 0 \equiv k \mod {32}
$$ 
with $\nabla Q(\x)$ not divisible by $4$. Once again, a lifting
argument would show that \eqref{N2} holds.
%T%REMOVED NEXT BIT
%T%To verify the claim about the solubility of the congruence, it clearly suffices to
%T%check all forms of the type
%T%$a_1 X_1^2 + \cdots + a_m X_m^2 + b X_{m+1}^2$,
%T%where $b=2$ or $b=6$ and $a_i \in \{1,3,5,7\} \; (1 \le i \le m)$.
We may therefore suppose that apart from $a_1, \ldots, a_m$, all of
 the coefficients of $Q$ are divisible by $4$. Let us suppose that
$m=4$, which is the worst case. Since $x_1,x_2,x_3, x_4$ 
are forced to be even in any $2$-adic solution, we may make the
 substitution $x_i=2x_i'$  for $1 \le i \le 4$. 
Moreover, we can carry out two reduction steps at this first stage, in
 the sense that after making these substitutions we get a new form
 $Q'$ which has all of its coefficients
divisible by $4$. Thus we may divide everything through by $4$. In
 this way we can replace the pair $(k,Q)$ by a pair $(k',Q')$, for which \eqref{31-pred}
holds with $p=2$ and $\theta=2$.  This iterative process 
clearly terminates after a finite number of steps, and produces a
pair $(k',Q')$ for which \eqref{31-pred} and \eqref{31-theta} holds, 
where the right-most term in \eqref{31-pred} can be bounded below as
$N(2^{t-\theta}) \gg 2^{(t-\theta)(n-1)}$.

On combining all of the findings above with \eqref{mendip1}, we have therefore
shown that
$$
\ss(k,Q) \gg_\ve|\D|^{-\ve} \prod_{p\mid \D}\frac{1}{p^{\nu_p(\D)/(n-4)}} \gg_\ve |\D|^{-1/(n-4)-\ve},
$$
when $n\geq 5$ and $(k,Q)$ satisfies weak \textsf{LSC}.
Alternatively, when $n\geq 5$ and $(k,Q)$ satisfies strong
\textsf{LSC}, we find that
$$
\ss(k,Q) \gg \prod_{p\mid \D}\Big(1-\frac{c_n}{p}\Big) \gg_\ve |\D|^{-\ve},
$$
for any $\ve>0$. These bounds are clearly satisfactory for Proposition
\ref{ss'}.

\end{document}